\documentclass{amsart}
\title[Tridiagonalized GUE matrices and labeled mobiles]{Tridiagonalized GUE matrices are a matrix model for labeled mobiles}
\author[A. Abdesselam, G. W. Anderson, and A. R. Miller]{Abdelmalek Abdesselam,  Greg W. Anderson$^*$ and Alexander R. Miller}
\address{University of Virginia Department of Mathematics, P. O. Box 400137, Charlottesville,
VA 22904-4137, USA}\email{malek@virginia.edu}
\address{School of Mathematics, University of Minnesota, Minneapolis, MN 55455, USA}\email{gwanders@umn.edu}
\thanks{$^*$Corresponding author}
\address{School of Mathematics, University of Minnesota, Minneapolis, MN 55455, USA}
\email{armiller@umn.edu}
\keywords{random matrices, tridiagonalized GUE matrices, BKAR formula, planar
maps, mobiles, cluster expansions, matrix models}
\subjclass[2010]{05C30, 05C10, 05E15, 15A52, 82B05}

\date{April 28, 2014}
\usepackage{graphicx,caption,float}
\usepackage[labelformat=simple]{subcaption}
\usepackage{mathrsfs}
\newcommand{\nn}{\langle n\rangle}

\newcommand{\mbold}{{\mathbf{m}}}
\newcommand{\LF}{{\mathrm{LF}}}
\newcommand{\GJdM}{{\mathrm{GJdM}}}
\newcommand{\dMotz}{{\mathrm{dMotz}}}
\newcommand{\Tri}{{\mathrm{Tri}}}
\newcommand{\GJ}{{\mathrm{GJ}}}

\newcommand{\norm}[1]{{\left\Vert#1\right\Vert}}
\newcommand{\Motz}{{\mathrm{Motz}}}
\newcommand{\Bond}{{\mathrm{Bond}}}
\newcommand{\Tree}{{\mathrm{Tree}}}
\def\bbone{{\mathchoice {\rm 1\mskip-4mu l} {\rm 1\mskip-4mu l}
{\rm 1\mskip-4.5mu l} {\rm 1\mskip-5mu l}}}
\newcommand{\indicator}[1]{\bbone\{#1\}}
\newcommand{\Map}{{\mathrm{Map}}}
\newcommand{\Orbit}{{\mathrm{Orb}}}
\newcommand{\supp}{{\mathrm{supp}}}
\newcommand{\zero}{{\mathbf{0}}}
\newcommand{\one}{{\mathbf{1}}}
\newcommand{\trace}{{\mathrm{tr}}}
\newcommand{\Cfrak}{{\mathfrak{C}}}
\newcommand{\Mat}{{\mathrm{Mat}}}
\newcommand{\Part}{{\mathrm{Part}}}

\newtheorem{Proposition}[subsubsection]{Proposition}
\newtheorem{Theorem}[subsubsection]{Theorem}
\newtheorem{Lemma}[subsubsection]{Lemma}
\newcommand{\RR}{{\mathbb{R}}}
\newcommand{\Sym}{{\mathrm{Sym}}}
\newcommand{\ZZ}{{\mathbb{Z}}}
\newcommand{\dd}{{\mathrm{d}}}
\newcommand{\Gfrak}{{\mathfrak{G}}}
\newcommand{\Tfrak}{{\mathfrak{T}}}
\newcommand{\Mfrak}{{\mathfrak{M}}}
\newcommand{\CC}{{\mathbb{C}}}
\newcommand{\Ebold}{{\mathbf{E}}}
\newcommand{\Qfrak}{{\mathfrak{Q}}}
\newcommand{\lf}{{\mathrm{LF}}}
\newcommand{\SV}{{\mathrm{SV}}}
\newcommand{\PP}{{\mathbb{P}}}
\begin{document}
\begin{abstract} 
It is well-known that the number of planar maps with prescribed vertex degree distribution
and suitable labeling
can be represented as the leading  coefficient
of the $\frac{1}{N}$-expansion of a joint cumulant of traces of powers of an 
  $N$-by-$N$ GUE matrix.
Here we undertake the calculation of this leading  coefficient in a different way.
Firstly,  we tridiagonalize the GUE matrix  {\em \`{a} la} \linebreak Trotter and Dumitriu-Edelman
and then alter it by conjugation to make the subdiagonal identically equal to $1$.
Secondly,  we apply the cluster expansion technique (specifically, the Brydges-Kennedy-Abdesselam-Rivasseau formula) from rigorous statistical mechanics.  Thirdly, by sorting through the terms of the expansion thus generated we arrive at an alternate interpretation for the leading coefficient related to factorizations of the long cycle $(12\cdots n)\in S_n$. Finally,
we reconcile the group-theoretical objects emerging from our calculation with the labeled mobiles of Bouttier-Di Francesco-Guitter.
\end{abstract}
\maketitle

\tableofcontents

\section{Introduction and main results}\label{section:MainResult}
Physicists in the 70's starting with 't Hooft \cite{tHooft} developed a beautiful combinatorial interpretation
for the limit
\begin{equation}\label{equation:TheLimit}
\lim_{N\rightarrow\infty} N^{\ell-2-\frac{n}{2}}\kappa\left(\trace\,\Xi_N^{\lambda_1},\dots,\trace\,\Xi_N^{\lambda_\ell}\right)
\end{equation}
where 
$\lambda$ is a partition, $n=|\lambda|$, $\ell=\ell(\lambda)$, 
$\Xi_N$ is an $N$-by-$N$ standard GUE matrix,
 and $\kappa(\cdot)$ is the joint cumulant functional.
 Namely, they interpreted \eqref{equation:TheLimit} as the number of suitably labeled planar maps with vertex degree distribution $\lambda$.
We recall details of this interpretation later in this introduction.
 
The goal of this paper is to recalculate the limit \eqref{equation:TheLimit} 
using a different toolbox to get a different interpretation for the same number.
The resulting interpretation 
counts objects related to factorizations of the long cycle $(12\cdots n)\in S_n$ which can then be put naturally in bijection with the labeled mobiles introduced by Bouttier-Di Francesco-Guitter \cite{BouFraGui} to enumerate
planar maps. We do not have anything new to say here about enumeration of planar maps
since the point of \cite{BouFraGui} is already to count them bijectively in terms of labeled mobiles.
Rather, the point is that our calculation of limit
\eqref{equation:TheLimit} using standard tools quite different from the usual ones leads naturally to the labeled mobiles---without any reference to planar maps whatsoever. In spirit (if not at all in the details) our work is similar to that of \cite{BouFraGui0} in that we count well-labeled trees
by approaching GUE matrices from a novel angle.

Of course the physicists went much farther and developed for all coefficients of the
 $\frac{1}{N}$-expansion of the joint cumulant appearing in the limit \eqref{equation:TheLimit} 
an interpretation in terms of higher genus maps. For simplicity we focus  in this paper exclusively
on the leading order, not venturing beyond genus zero.  However, our method is not intrinsically
limited to genus zero.

The prevailing view in combinatorics is to favor bijective proofs over other less constructive ways
of establishing that two finite sets have the same cardinality.
From that point of view our  achievement is less than stunning: in effect, we give a long analytic proof
for a fact which already has been given a fairly short bijective proof in \cite{BouFraGui}.  But our method of proof yields a connection between labeled mobiles and tridiagonalized GUE matrices of intrinsic interest which might possibly serve  as a heuristic device in ways yet to be worked out. For example, the graph metric playing such an important role in \cite{BouFraGui} appears naturally in the tridiagonal context when one ``opens the brackets'' as in \S\ref{subsubsection:OpenBrackets} below;
possibly a hierarchy of labelings relevant to enumeration of higher genus maps could be discovered by continuing the analysis started here.

Here is an outline of our recomputation of \eqref{equation:TheLimit}.
\begin{enumerate}
\item[(I)]
We replace the $N$-by-$N$ GUE matrix $\Xi_N$ by its (lightly modified) tridiagonalization 
 {\em \`{a} la} Trotter \cite{Trotter} and Dumitriu-Edelman \cite{DumitriuEdelman1}. We carry out this  easy step in \S\ref{subsection:ThePlan} below.
See \cite[Section 4.5]{AGZ} for background on tridiagonalization.
 \item[(II)] We apply the Brydges-Kennedy-Abdesselam-Rivasseau (BKAR) formula
from rigorous statistical mechanics to obtain a delicate expansion 
of the joint cumulant under the limit in \eqref{equation:TheLimit}. We carry out this step in 
\S\ref{subsection:FieldObservation} below after setting up the BKAR machinery in the preceding part of \S\ref{section:BKAR}. The formula in question is \eqref{equation:MotzkinExpansionSteroid} below.
\item[(III)] We analyze the many terms summed up in formula \eqref{equation:MotzkinExpansionSteroid} in order to obtain our main result,  which is phrased in group-theoretical terms.
The calculations in question
are carried out in \S\ref{section:GraphPreparation} and \S\ref{section:TridiagonalCumulantCalc} below.
Our main result appears as Theorem \ref{Theorem:MainResult} below.
\end{enumerate}
 It is worth remarking that this is probably the first paper in which the BKAR formula has been used to perform
an exact calculation; ordinarily one applies it only to obtain upper bounds on joint cumulants.
Since we provide background and references for the BKAR formula in \S\ref{section:BKAR},
along with a short proof, we omit  further discussion of it in this introduction.

Our main result interprets the limit \eqref{equation:TheLimit} as the cardinality of a certain group-theoretically defined  set. The objects so counted have the following notable features:
\begin{itemize}

\item They may be identified with the labeled mobiles of  Bouttier-Di Francesco-Guitter
\cite{BouFraGui}.  (See  \S\ref{subsection:PMM} below for details.)
\item They thus belong to a line of research developed over many decades focused on bijections
between sets of well-labeled trees and planar maps. See, e.g.,  \cite{AmbjornBudd},  \cite{BerFus}, 
\cite{BousJeh}, \cite{BousScha}, \cite{BouFraGui}, \cite{CoriVauquelin}, \cite{Schaeffer}
and \cite{Tutte}.
A driver of research into such bijections lately has been the intense activity in probability and physics in connection with the  Brownian map and quantum gravity.
Concerning the latter see, e.g., \cite{BetJacMier}, \cite{ChaSch},
\cite{LeGall}, \cite{Miermont}.
\item  They can be analyzed with the help of a combinatorial insight of Goulden-Jackson
\cite{GouldenJackson} concerning factorizations  of the long cycle $(12\cdots n)\in S_n$.
See \S\ref{subsection:TutteRecovery} below for this analysis.
\item They are the group-theoretical counterparts of simple examples 
of \linebreak Grothendieck's {\em dessins d'enfants} \cite{SchnepsEtAl}, \cite{ShabatVoevodsky},
albeit with extra ``coloring.''
\item They make it possible to straightforwardly reconcile
 our  interpretation of \eqref{equation:TheLimit}
with a famous formula of Tutte \cite{TutteSlice} for the number of Eulerian (all degrees even)
rooted planar maps with prescribed vertex degree distribution. 
See \S\ref{subsection:PlanarMapIntro} and \S\ref{subsection:TutteRecovery} below for details.
\end{itemize}
Of course the last point is hardly surprising given the results of \cite{BouFraGui}.
We work out the exercise of recovering Tutte's formula in order to warm the reader up
for the comparison with the theory of \cite{BouFraGui} undertaken in \S\ref{subsection:PMM}.

In the remainder of this (rather long) introduction we formulate our main result precisely
and provide  details concerning several points briefly mentioned above.

\subsection{Table of notation}\label{subsection:GeneralNotation}
We briefly mention the most basic items of notation and terminology used throughout the paper.
The reader should scan the table once quickly and then  use it as a reference.

\subsubsection{General notation and terminology}
 Let $|S|$ denote the cardinality of a finite set $S$.
  Let $\indicator{\cdot}$ be probabilist's indicator notation.
The $(i,j)$-entry of a matrix $A$ is invariably denoted by $A(i,j)$.
Let $\langle n\rangle=\{1,\dots,n\}$ for positive integers $n$.
 Let $\Part_n$ denote the lattice of partitions of the set $\langle n\rangle$. (For further notation related to $\Part_n$, see \S\ref{section:BKAR} below.)
Constants in estimates are usually denoted by $c$, $C$, or $K$,
and their numerical values may change from line to line.

\subsubsection{Numerical partitions}
A {\em numerical partition} (or simply {\em partition}, context permitting) is a monotone decreasing sequence $\lambda=\{\lambda_i\}_{i=1}^\infty$ of nonnegative integers
such that $\lambda_i=0$ for $i\gg 0$. The (nonzero) terms $\lambda_i$ are called 
the {\em parts} of $\lambda$.  By and large we follow notation of Macdonald \cite{Macdonald}.  Let $|\lambda|=\sum_i \lambda_i$. 
We also write $\lambda\vdash n\Leftrightarrow |\lambda|=n$. 
 Let $m_i(\lambda)=|\{j\mid \lambda_j=i\}|$ for $i>0$.
 Let $\ell(\lambda)=\sum_{i}m_i(\lambda)=|\{i\mid \lambda_i>0\}|$, called the {\em length} of $\lambda$.
 Let $z_\lambda=\prod_ii^{m_i(\lambda)}m_i(\lambda)!$. Abusing notation we occasionally write
 $\lambda=\prod_i i^{m_i(\lambda)}$.

\subsubsection{Graphs}
For us a {\em graph} is a finite set of vertices and a finite set of edges along with
the specification of an incidence relation which designates for each edge a
set of one or two endpoints among the vertices.  Furthermore, in the case of graphs
without multiple edges, we simply identify edges with their endpoint sets.

\subsubsection{Permutations}
Let $S_n$ denote the group of permutations of $\langle n\rangle$. 
 For $\sigma\in S_n$, let 
$\supp\, \sigma=\{i\in \langle n\rangle\mid \sigma(i)\neq i\}$, which we call the {\em support} of $\sigma$.
In other words, $\supp\,\sigma$ is the complement of the set of fixed points of $\sigma$.
A {\em cycle} in $S_n$ is a permutation with nonempty support on which it acts transitively. The {\em length} of a cycle is the cardinality of its support.  Cycles are called {\em disjoint} if they have disjoint supports.
A cycle of length $m$ is called an {\em $m$-cycle}.  
(In our usage there are no $1$-cycles.)
A $2$-cycle is also called a {\em transposition}.  
  Each $\sigma\in S_n$ has a 
 factorization into disjoint cycles unique up to ordering of the factors, hereafter called simply the {\em canonical factorization} of $\sigma$.
For $\sigma\in S_n$, let $\Orbit_n(\sigma)\in \Part_n$ denote the finest partition
consisting of \linebreak $\sigma$-stable blocks. Blocks of $\Orbit_n(\sigma)$ are called {\em $\sigma$-orbits}.
(Whereas $1$-cycles are disallowed here, $\sigma$-orbits may of course be singletons.)
As usual we index the conjugacy class of a permutation $\sigma\in S_n$ by the numerical partition $\lambda\vdash n$ 
recording  the cardinalities of the blocks of the set partition $\Orbit_n(\sigma)$.
We also write  $\ell(\sigma)=|\Orbit_n(\sigma)|=\ell(\lambda)$ and $\sigma\sim \lambda$.

\subsection{A statement of the main result}
We now formulate the main result of the paper in terms of permutations only.
This way of presenting the result is brief but it is also  misleading, as we will explain presently.

\subsubsection{The set $\Map_n$}
Let 
$\Map_n$ denote the set of ordered pairs 
$(\theta,\iota)\in S_n\times S_n$
of permutations satisfying the following conditions:
\begin{eqnarray}\label{equation:PM1}
&&\mbox{$\iota$ is fixed-point-free and squares to the identity.}\\
\label{equation:PM3}
&&\ell(\theta)-\ell(\iota)+\ell(\theta\iota)=2, \;\mbox{cf. Euler's formula $V-E+F=2$.}\\
\label{equation:PM4}
&&\mbox{$\theta$ and $\iota$ generate a subgroup of $S_n$ acting transitively on $\langle n\rangle$.}
\end{eqnarray}
For convenience we also define
$$\Map_n(\theta)=\{\iota\in S_n\mid (\theta,\iota)\in \Map_n\}.
$$
Clearly we have
\begin{equation}\label{equation:AintZero}
\frac{n}{2}-\ell(\theta)+2\leq 0\Rightarrow \Map_n(\theta)=\emptyset.
\end{equation}
The set $\Map_n$ is allied with planar maps in a fashion we recall
briefly in \S\ref{subsection:PlanarMapIntro} below.

\subsubsection{The set $\GJ_n$} 
Let $\GJ_n$ denote the set of ordered pairs 
$(\theta,\sigma)\in S_n\times S_n$
of permutations satisfying
\begin{equation}
\label{equation:GouldenJackson}
\ell(\theta)+\ell(\sigma)=n+1\;\mbox{and}\;
\ell(\theta\sigma)=1.
\end{equation}
Members of $\GJ_n$ will be called {\em Goulden-Jackson} pairs. 
For convenience we define
$$\GJ_n(\theta)=\{\sigma\in S_n\mid (\theta,\sigma)\in \GJ_n\}.$$
In \S\ref{subsection:TutteRecovery} we recall the  interpretation of elements of
$\GJ_n$ in terms of planar trees.

\subsubsection{The set $\dMotz_n(\theta,\sigma)$}
Given $(\theta,\sigma)\in \GJ_n$, 
let $\dMotz_n(\theta,\sigma)$ denote the set of functions 
$g:\langle n\rangle\rightarrow\ZZ$
with the following properties:
\begin{eqnarray}
\label{equation:Interaction0}
&&|g|\leq 1.\\
\label{equation:Interaction1}
&& \mbox{$g$ averages to $0$ on $\theta$-orbits.}\\
\label{equation:Interaction2}
&&g\circ \sigma=g.\\
\label{equation:Interaction3}
&&\{g=-1\}\cap\supp\,\sigma=\emptyset.\\
\label{equation:Interaction3.5}
&&\{g=0\}\cap \supp\,\sigma^2=\emptyset.\\
\label{equation:Interaction4}
&&\{g=0\}\subset \supp\,\sigma.
\end{eqnarray}
The rationale  for the  (ungainly) notation $\dMotz$
is given in Proposition \ref{Proposition:EquiNumerous} below.
\subsubsection{The sets $\GJdM_n$ and $\GJdM_n(\theta)$}
Combining notions introduced above, we define the following more complicated sets:
\begin{eqnarray}
\GJdM_n
&=&\{(\theta,\sigma,g)\in \GJ_n\times \{0,\pm 1\}^{\langle n\rangle}\mid 
g\in \dMotz_n(\theta,\sigma)\}.\\
\GJdM_n(\theta)
&=&\{(\sigma,g)\in S_n\times \{0,\pm 1\}^{\langle n\rangle}\mid (\theta,\sigma,g)\in \GJdM_n\}.
\end{eqnarray}
 We will later show that
\begin{equation}\label{equation:AintZeroNeither}
\frac{n}{2}-\ell(\theta)+2\leq 0\Rightarrow \GJdM_n(\theta)= \emptyset.
\end{equation}
(See Lemma \ref{Lemma:BalanceBooksBis} below.)
We briefly indicate in \S\ref{subsection:TutteRecovery} below a graphical interpretation 
for members of $\GJdM_n$ in terms of vertex-four-colored edge-labeled planar trees.
Furthermore, in \S\ref{subsection:PMM} we will explain how to identify these objects with the
labeled mobiles of \cite{BouFraGui}.

Here is the main result of the paper.

\begin{Theorem}\label{Theorem:MainResult}
For all $\theta\in S_n$ such that $\frac{n}{2}-\ell(\theta)+2>0$ we have
\begin{equation}\label{equation:MainResult}
|\Map_n(\theta)|=\frac{|\GJdM_n(\theta)|}{\frac{n}{2}-\ell(\theta)+2}.
\end{equation}
\end{Theorem}
\noindent By relations \eqref{equation:AintZero} and \eqref{equation:AintZeroNeither} noted above, the
 numerical hypothesis $\frac{n}{2}-\ell(\theta)+2>0$ is merely a convenience excusing us from having to break out trivial cases for separate examination, not an essential restriction.
 The proof of Theorem \ref{Theorem:MainResult}
commences in \S\ref{section:BKAR} and takes up the rest of the paper.
 \subsubsection{Remark}
 Theorem \ref{Theorem:MainResult} is indeed misleading in the simplified combinatorial presentation given above.
First of all, in light of the results of \cite{BouFraGui} and the possibility of interpreting elements of $\GJdM_n$
as labeled mobiles discussed in \S\ref{subsection:PMM} below, Theorem \ref{Theorem:MainResult} as stated is not new---we claim novelty only for our proof, which is analytic, proceeding by way of the study of tridiagonalized GUE matrices. See Remark \S\ref{subsubsection:ThePoint} below for further discussion of this point. 
Secondly, Theorem \ref{Theorem:MainResult} conceals what we consider to be the most important contribution of the paper.
Namely, we feel considerable value attaches to the possibility of extending our calculations to higher genus, leading perhaps to the discovery (if not the proof) of new ways of enumerating higher genus maps. See Proposition \ref{Proposition:ProgressSoFar} below for an exact formula
which we believe has a $\frac{1}{N}$-expansion worth working out to all orders.

\subsection{Planar maps and Tutte's formula}
\label{subsection:PlanarMapIntro}
We recall intuitions guiding the study of the set $\Map_n$,
introduce notation needed throughout the paper
and finally recall a famous result of Tutte.

\subsubsection{The link between planar maps and permutation pairs}
A {\em planar map} is a cellular decomposition of the $2$-sphere with connected $1$-skeleton.
We call $0$-cells (resp., $1$-cells and $2$-cells) {\em vertices} (resp., {\em edges} and {\em faces}).
Each edge is viewed as two half-edges stuck together.
The {\em degree} of a vertex is the number of half-edges incident upon it.
A {\em half-edge-labeled} planar map of $n$ half-edges is a planar map 
equipped with a numbering from $1$ to $n$ of its half-edges.
Each half-edge-labeled planar map of $n$ half-edges gives rise to a permutation pair $(\theta,\iota)\in \Map_n$
by the following procedure.  Let $\theta$ be the permutation whose cycles record in counterclockwise order the labels of half-edges sprouting from the vertices of degree $>1$; every label of a half-edge terminating in a vertex of degree $1$ is a fixed point of $\theta$. Let $\iota$ be the permutation which exchanges labels of half-edges
belonging to the same edge. It is well-known that every element $(\theta,\iota)\in\Map_n$ arises from
a half-edge-labeled planar map in the manner just specified. We regard two half-edge-labeled planar maps as equivalent if both give rise to the same element of $\Map_n$.
See the series of three survey papers \cite{CoriMachi}
for an introduction to the  point of view emphasizing permutation pairs.
Generally our attitude is that permutation pairs are the objects of rigorous study in this paper,
whereas we view planar maps and related graphs as 
(very appealing) heuristic devices.

\begin{figure}[hbt]
\centering
\includegraphics[scale=1]{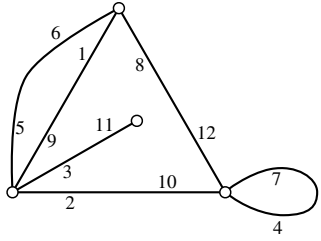}
\caption{
This drawing represents a half-edge-labeled planar map corresponding to 
the pair $(\theta,\iota)\in\Map_{12}$ where $\theta=(1,8,6)(2,3,9,5)(4,7,12,10)$ and $\iota=(1,9)(2,10)(3,11)(4,7)(5,6)(8,12)$.  The vertex degree distribution of this planar map  is the numerical partition $4^2\cdot 3^1\cdot 1^1$.
\label{HalfEdges:Fig}
}
\end{figure}

\subsubsection{Rooted planar maps}
A {\em rooted planar map} of $n$ half-edges is
(in effect) a half-edge-labeled planar map from which one erases all of the labels but $n$.
Let us identify $S_{n-1}$ with the subgroup of $S_n$ consisting of permutations
fixing the point $n$ and let $S_{n-1}$ act on $\Map_n$ by simultaneous conjugation,
i.e., the action of $\rho\in S_{n-1}$ on $(\theta,\iota)\in \Map_n$ is $(\rho\theta\rho^{-1},\rho\iota\rho^{-1})\in \Map_n$. 
Then rooted planar maps (up to equivalence) are indexed by the orbit space $\Map_n/S_{n-1}$.

\begin{figure}[hbt]
\centering
\includegraphics[scale=1]{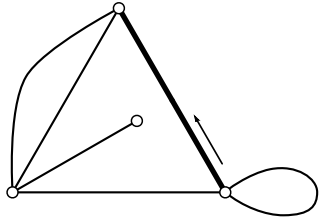}
\caption{
This drawing shows the rooted planar map arising by erasure of labels from the half-edge-labeled planar map depicted in Figure 
\ref{HalfEdges:Fig}.
\label{RootedPlanarMap:Fig}
}
\end{figure}

\begin{Lemma}\label{Lemma:NoAutomorphisms} $S_{n-1}$ acts freely on $\Map_n$.
\end{Lemma}
\proof  Fix $\rho\in S_{n-1}$ and $(\theta,\iota)\in \Map_n$ such that
$(\theta,\iota)=(\rho\theta\rho^{-1},\rho\iota\rho^{-1})$, i.e., such that $\rho$ commutes with both $\theta$ and $\iota$. It is enough to show that $\rho=1$. In any case,  the set of points of $\langle n\rangle$ fixed by $\rho$ is not empty
and moreover stable under the action of the group of permutations generated by $\theta$ and $\iota$.
But the latter group by definition of $\Map_n$ acts transitively on $\langle n\rangle$.
Thus every point of $\langle n\rangle$ is fixed by $\rho$.
\qed

\subsubsection{The numbers $\Mfrak_\lambda$ and $\Mfrak_\lambda^\star$}
Let $\lambda\vdash n$ be a partition and let $\ell=\ell(\lambda)$.
For any permutation $\theta\in S_n$ belonging
to the conjugacy class indexed by $\lambda$ let 
$$\Mfrak_\lambda=|\Map_n(\theta)|.$$
The number $\Mfrak_\lambda$ is well-defined because the number on the right depends only on the conjugacy class of $\theta$.
Let $\Map_{\lambda\vdash n}$ denote the subset of $\Map_n$ consisting of pairs $(\theta,\iota)$
such that $\theta$ belongs to the conjugacy class indexed by $\lambda$. 
Clearly $\Map_{\lambda\vdash n}$ is stable under the action of $S_{n-1}$.
Let
$$\Mfrak_\lambda^\star=|\Map_{\lambda\vdash n}/S_{n-1}|.$$
The number $\Mfrak_\lambda^\star$ counts (equivalence classes of) rooted planar maps
having  vertex degree distribution $\lambda$.
Since $\frac{n!}{z_\lambda}$ is the cardinality of the conjugacy class in $S_n$ indexed by $\lambda$
and $S_{n-1}$ acts freely on $\Map_n$ by Lemma \ref{Lemma:NoAutomorphisms},
one has the comparison formula
\begin{equation}\label{equation:StarConversion}
\Mfrak_\lambda=\frac{z_\lambda}{n}\Mfrak^\star_\lambda.
\end{equation}
Thus the numbers $\Mfrak_\lambda$ and $\Mfrak^\star_\lambda$ 
 carry the same information even if they have rather different connotations.

\subsubsection{Counting of Eulerian rooted planar maps, following Tutte}
Recall that planar maps with all vertex degrees even are called {\em Eulerian}.
Now let $\lambda\vdash n$ and $\ell=\ell(\lambda)$, as above. Also let $m_i=m_i(\lambda)$.
Assume that every part $\lambda_i$ is even so that $\lambda$ is a possible vertex-degree distribution
of an Eulerian planar map.
Tutte \cite{TutteSlice}  has given in the Eulerian case a simple explicit formula for the number $\Mfrak^\star_\lambda$, namely
\begin{equation}\label{equation:Tutte}
\Mfrak_\lambda^\star
=
\frac{2(\frac{n}{2})!}{(\frac{n}{2}-\ell+2)!}
\cdot \prod_{i\geq 1} \frac{1}{m_{2i}!}\left(\begin{array}{c}
2i-1\\
i\end{array}\right)^{m_{2i}}.
\end{equation}
See also \cite{Schaeffer} for a more recent proof of this same formula by an elegant
construction of a bijection.
Using \eqref{equation:StarConversion} above we can rewrite Tutte's formula \eqref{equation:Tutte} as
\begin{equation}\label{equation:TutteBis}
\Mfrak_\lambda=\frac{(\frac{n}{2}-1)!}{(\frac{n}{2}-\ell+2)!}\cdot \prod_{i=1}^\ell \frac{\lambda_i}{2}
\cdot
\prod_{i=1}^\ell\left(\begin{array}{c}
\lambda_i\\
\lambda_i/2\end{array}\right).
\end{equation}
We will find the latter presentation of Tutte's result more convenient.

\subsection{Enumeration of planar maps via matrix integrals}\label{subsection:Physicists}
We turn next to the physicists' point of view on the numbers $\Mfrak_\lambda$.

\subsubsection{Standard GUE matrices}
A random $N$-by-$N$ hermitian matrix $\Xi$ is called a {\em standard GUE matrix}
if its law has   the density 
$\exp\left(-\frac{1}{2}\trace \,H^2\right)$ with respect to Lebesgue measure, up to a normalization factor.
Equivalently, one requires the family $\{\Xi(i,j)\}_{1\leq i\leq j\leq N}$ of matrix entries on or above the diagonal to be independent and to have a centered Gaussian joint
distribution characterized by $\Ebold \Xi(i,j)^2=\delta_{ij}$
and $\Ebold |\Xi(i,j)|^2=1$.

\subsubsection{The number $\Mfrak_{\lambda,N}$ and its leading order behavior}
Let $\lambda$ be a numerical partition and let $\ell=\ell(\lambda)$. 
Let $N$ be a positive integer. Let $\Xi_N$ be a standard $N$-by-$N$ GUE matrix.
 Let
\begin{equation}\label{equation:FiniteMfrakN}
\Mfrak_{\lambda,N}= 
\kappa\left(\trace\, \Xi_N^{\lambda_1},\dots,\trace\, \Xi_N^{\lambda_\ell}\right)
\end{equation}
where $\kappa(\cdot)$ is the joint cumulant functional. (See \S\ref{subsection:JointCumulants} 
below to be reminded of the definition and first properties of joint cumulants.)
Physicists in the 1970's obtained the limit formula
\begin{equation}\label{equation:tHooft}
\Mfrak_\lambda=\lim_{N\rightarrow\infty}
N^{\ell-2-\frac{n}{2}}\Mfrak_{\lambda,N}.
\end{equation}
The right side here is precisely the limit \eqref{equation:TheLimit} with which we began the introduction.
More generally physicists derived for $\Mfrak_{\lambda,N}$ an asymptotic expansion in
powers of $\frac{1}{N}$ with coefficients counting diagrams of higher genus.
But in this paper we will be content to study genus zero (leading order) behavior only. 

\begin{figure}[hbt]
\centering
\includegraphics[scale=0.8]{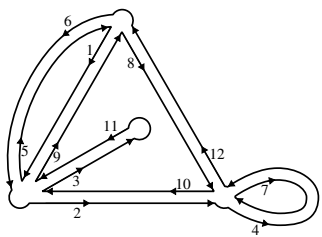}
\caption{This drawing is a rendering of Figure 1 as a {\em fat graph}
or {\em ribbon graph}.
\label{Fig:FatGraphS}
}
\end{figure}

\subsubsection{Notes and references}
The paper \cite{tHooft} is recognized as the initiation of GUE enumeration of maps
 although no formula recognizable to a mathematician as \eqref{equation:tHooft} could be found there. 
Many hands subsequently developed the theory around formula \eqref{equation:tHooft}.
Without any pretension to completeness,
we mention the references \cite{BIPZ}, \cite{ErcolaniMcLaughlin}, \cite{GuionnetMaurelSegala},
\cite{LandoZvonkin} and \cite{Zvonkin}  as ways to enter this vast territory.
\subsubsection{Remark}
Using Weisner's theorem \cite[p. 351]{Rota}
in conjunction with the Isserlis-Wick formula \eqref{equation:WickFormula} recalled below,
it is possible to give a proof of  formula \eqref{equation:tHooft} 
completely within the domain of algebraic combinatorics,
using the permutation pair point of view.

\subsubsection{Remark}\label{subsubsection:ThePoint}
To prove Theorem \ref{Theorem:MainResult} we will deal with the number $\Mfrak_\lambda$
solely through formula \eqref{equation:tHooft}.
Without loss of comprehension, from \S\ref{section:BKAR} of the paper onward, the reader could take formula \eqref{equation:tHooft} as the definition of $\Mfrak_\lambda$.  The point of the paper really is to provide a re-interpretation of the limit ``from scratch,'' including a self-contained proof of the existence of the limit.
Planar maps will be nowhere in sight.

\subsection{The tridiagonal representation of $\Mfrak_{\lambda,N}$}
\label{subsection:ThePlan}
We carry out the first step of our recalculation of the limit \eqref{equation:TheLimit}.

\subsubsection{Tridiagonalization of standard GUE matrices}
Let $N$ be a positive integer which eventually we send to infinity.
Let $\Xi_N$ be an $N$-by-$N$ standard GUE matrix.
The result of applying to $\Xi_N$ the well-known Givens-Householder tridiagonalization procedure (albeit starting at the lower right corner
rather than the upper left) yields a random matrix with the same law
as that of the random matrix
\begin{equation}\label{equation:Depicted}
\left[\begin{array}{cccccccccc}
\xi_1&\sqrt{\eta_1}\\
\sqrt{\eta_1}&\xi_2&\sqrt{\eta_2}\\
&\sqrt{\eta_2}&\xi_3&\sqrt{\eta_3}\\
&&\sqrt{\eta_3}&\ddots&\ddots\\
&&& \ddots   &\xi_{N-1}&\sqrt{\eta_{N-1}}\\
&&&&\sqrt{\eta_{N-1}}&\xi_{N}
\end{array}\right],
\end{equation}
where the family
$$
\{\xi_i\}_{i=1}^N\cup\{\eta_i\}_{i=1}^{N-1}
$$
of real random variables is independent, each random variable $\xi_i$ is standard normal
and each random variable $\eta_i$ has a 
distribution of $\Gamma$-type for which
$$\frac{e^{-x}x^{i-1}}{(i-1)!}\indicator{x>0}$$
is the probability density function. 
The tridiagonalization procedure does not change the law of the spectrum,
and so in the definition \eqref{equation:FiniteMfrakN} 
one could in principle replace $\Xi_N$ 
by the random matrix \eqref{equation:Depicted}.

\subsubsection{Notes and references on tridiagonalization}
The idea to approach the semicircle law for GOE matrices through tridiagonalization
is due to Trotter \cite{Trotter}. The idea of tridiagonalization was later developed to yield
tridiagonal matrix models for Gaussian $\beta$-ensembles by Dumitriu-Edelman \cite{DumitriuEdelman1}
and developed further by these authors to yield a CLT \cite{DumitriuEdelman2}. 
The matrix \eqref{equation:Depicted} is the case $\beta=2$ of the Dumitriu-Edelman model
for the Gaussian $\beta$-ensemble.
Here, for simplicity,  we do not push our analysis beyond the case of $\beta=2$. 
See \cite[Section 4.5]{AGZ} for background on tridiagonalization and $\beta$-ensembles.

\subsubsection{Another tridiagonal representation of $\Mfrak_{\lambda,N}$}
Now the random variables $\sqrt{\eta_i}$ are not so nice for our purposes.
Accordingly, we conjugate
the matrix \eqref{equation:Depicted} suitably
to get another tridiagonal matrix
\begin{equation}\label{equation:DepictedTri}
\Tri_N=\left[\begin{array}{cccccccccc}
\xi_1&\eta_1\\
1&\xi_2&\eta_2\\
&1&\xi_3&\eta_3\\
&&1&\ddots&\ddots\\
&&& \ddots   &\xi_{N-1}&\eta_{N-1}\\
&&&&1&\xi_{N}
\end{array}\right]
\end{equation}
the spectrum of which has the same law. 
Replacing $\Xi_N$ in formula \eqref{equation:FiniteMfrakN} by $\Tri_N$ we obtain the formula
\begin{equation}\label{equation:PreMotzkinExpansion}
\Mfrak_{\lambda,N}=\kappa\left(\trace\, \Tri_N^{\lambda_1},\dots,\trace\,\Tri_N^{\lambda_\ell}\right)
\end{equation}
holding for every partition $\lambda$, where as usual $\ell=\ell(\lambda)$. 

\subsubsection{Further details on the plan of proof}
To prove Theorem \ref{Theorem:MainResult}
we will use formula \eqref{equation:PreMotzkinExpansion}
rather than formula \eqref{equation:FiniteMfrakN}
to evaluate the limit  on the right side of \eqref{equation:tHooft}.
The perhaps unexpected extra ingredient in our calculation is the BKAR formula from rigorous statistical mechanics.
(See Theorem \ref{Theorem:BKARarboreal} and its application Theorem \ref{Theorem:MainTool} below.) 
The BKAR formula will permit us to control cancellation on the right side of formula \eqref{equation:PreMotzkinExpansion} by means of repeated integration by parts.
The evaluation of the limit on the right side of \eqref{equation:tHooft} by the tridiagonal/BKAR route will take up
the rest of the paper from \S\ref{section:BKAR} onward.

\subsection{Recovery of Tutte's formula from Theorem \ref{Theorem:MainResult}}\label{subsection:TutteRecovery}
Our goal here is to reconcile formulas \eqref{equation:MainResult} and \eqref{equation:TutteBis}
by showing directly that their right sides are equal in the Eulerian case. The calculations needed to do this are completed in \S\ref{subsubsection:Snip} below after suitable preparation. 
In contrast to our mostly analytical {\em modus operandi} in this paper, under this heading and the next we use a relatively  informal ``tree-surgical'' approach.
As a byproduct of our discussion under this heading 
we provide a simple graphical interpretation for each member of the set $\GJdM_n$
which subsequently in \S\ref{subsection:PMM} below we explain how to view as a labeled mobile.

\subsubsection{Shabat-Voevodsky trees} The simplest examples
of the {\em dessins d'enfants} introduced by Grothendieck 
(see \cite{SchnepsEtAl} for background) are the two-colored planar trees.
These objects come to number-theoretic life in connection with the Shabat-Voevodsky polynomials \cite{ShabatVoevodsky}.
See also the short note \cite{Biane} 
for a simple and beautiful if not entirely elementary
construction of these polynomials.
 We will not delve into the theory of {\em dessins d'enfants} here, but we will acknowledge the tangential relationship of our work to this theory
by calling a bipartite edge-labeled planar tree (vertices colored white and black, with no two adjacent vertices of the same color,  and with edges numbered from $1$ to $n$, where $n$ is the number of edges)
a {\em Shabat-Voevodsky tree}.
Let $\SV_n$ denote the set of (equivalence classes of) Shabat-Voevodsky trees of $n$ edges. 

\subsubsection{Generalized definitions}
For technical flexibility we need to generalize several definitions
given above in a harmless way. Let $A$ be any finite set and let $n=|A|$.
Let $S_A$ denote the group of permutations of the set $A$.
Let $\SV_A$ denote the set of (equivalence classes
of) Shabat-Voevodsky trees with $n$ edges labeled by 
distinct elements of the set $A$ rather than by distinct elements of the set $\langle n\rangle$. In the same spirit, let $\GJ_A\subset S_A\times S_A$ 
denote the subset defined by evident analogy with the definition
of $\GJ_n$ in the case  $A=\langle n\rangle$.

\begin{figure}[tbh]
\centering
\begin{subfigure}[t]{0.32\textwidth}
\centering\includegraphics{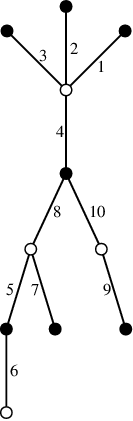}
\caption{}
\label{SVTree:Fig}
\end{subfigure}
\begin{subfigure}[t]{0.32\textwidth}
\centering\includegraphics{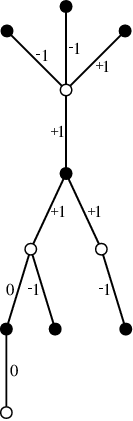}
\caption{}
\label{DJdM:Fig}
\end{subfigure}
\begin{subfigure}[t]{0.32\textwidth}
\centering\includegraphics{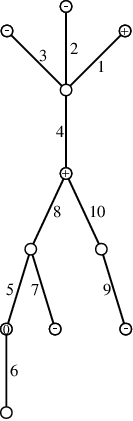}
\caption{}
\label{DJdMUp:Fig}
\end{subfigure}
\caption{
Let $\theta=(1,2,3,4)(5,7,8)(9,10)$ and $\sigma=(4,8,10)(5,6)$.  
Then $(\theta,\sigma)\in \mathrm{GJ}_{10}$.  
The corresponding Shabat-Voevodsky tree in $\mathrm{SV}_{10}$ is~\subref{SVTree:Fig}, 
and~\subref{DJdM:Fig} illustrates an element $g$ from $\mathrm{dMotz}_{10}(\theta,\sigma)$.  
The triple 
$(\theta,\sigma;g)\in \mathrm{GJdM}_{10}$
is encoded by $\mathfrak T_{10}(\theta,\sigma;g)$ in~\subref{DJdMUp:Fig}.
}\label{Diagrams:Fig}
\end{figure}

\begin{Proposition}\label{Proposition:GJSV} For finite sets $A$,
the sets $\SV_A$ and $\GJ_A$ are canonically in bijection.
\end{Proposition}
\noindent This is a commonplace both in the theory of {\em dessins d'enfants}
and in combinatorics in relation to the problem of calculating connection coefficients for conjugacy classes of the symmetric group.
\proof Given a Shabat-Voevodsky tree $\Tfrak$ belonging to $\SV_A$,
by writing down for each white vertex of degree $>1$ in counterclockwise order the labels of edges incident on the vertex, one obtains a permutation $\theta_\Tfrak\in S_A$ canonically factored into cycles;
the label of each edge terminating in a white leaf
is a fixed point of $\theta_\Tfrak$. Similarly one obtains a permutation $\sigma_\Tfrak\in S_A$
by reversing the roles of white and black.
One checks immediately that $(\theta_\Tfrak,\sigma_\Tfrak)\in \GJ_A$,
and that every $(\theta,\sigma)\in \GJ_A$ so arises in an essentially
unique way. See drawing (a) in Figure \ref{Diagrams:Fig} below for an illustration of the passage from a
Goulden-Jackson pair to a Shabat-Voevodsky tree.
The notation $\theta_\Tfrak$ and $\sigma_\Tfrak$ introduced in this proof
will be needed below to complete the job of reconciling Theorem \ref{Theorem:MainResult} with Tutte's formula.  
\qed

\begin{Lemma}\label{Lemma:GJEnumeration}
For partitions $\lambda\vdash n$ and $\theta\in S_n$ such that $\theta\sim\lambda$
we have 
\begin{equation}
|\GJ_n(\theta)|=\frac{(n-1)!}{(n-\ell+1)!}\cdot \prod_{i=1}^\ell \lambda_i
\end{equation}
where $\ell=\ell(\lambda)$.
\end{Lemma}
\proof Let $\mu$ be any partition such that
$\ell(\lambda)+\ell(\mu)=n+1$ and $|\lambda|=|\mu|=n$.
The result \cite[Thm. 2.2]{GouldenJackson} translated into the present setup
says that
$$
|\{(\rho,\sigma)\in S_n\times S_n\mid \rho\sigma=(1\cdots n),\;\rho\sim \lambda\;\mbox{and}\;\sigma\sim \mu\}|=
n\frac{(\ell(\lambda)-1)!(\ell(\mu)-1)!}{\prod_i m_i(\lambda)!\prod_j m_j(\mu)!}.
$$
We note that the statement above is originally due to other authors (see \cite{BedardGoupil})
and that it was originally proved by an inductive method.
We note also that the main goal of \cite{GouldenJackson}  was to give a different bijective proof 
of the same result. The  idea animating the latter proof
we have recapitulated as Proposition \ref{Proposition:GJSV} above.
It follows that
\begin{eqnarray*}
&&|\{(\rho,\sigma)\in \GJ_n\mid \rho\sim \lambda\;\mbox{and}\;\sigma\sim \mu\}|\;=\;
n!\frac{(\ell(\lambda)-1)!(\ell(\mu)-1)!}{\prod_i m_i(\lambda)!\prod_j m_j(\mu)!},\;\;\mbox{hence}\\
&&|\{\sigma\in \GJ_n(\theta)\mid \sigma\sim \mu\}|\;=\;z_\lambda\cdot\frac{(\ell(\lambda)-1)!(\ell(\mu)-1)!}{\prod_i m_i(\lambda)!\prod_j m_j(\mu)!}\;\;\mbox{and finally}\\
&&|\GJ_n(\theta)|
=(\ell(\lambda)-1)!\cdot (n-\ell(\lambda))!\cdot \prod_{i=1}^{\ell(\lambda)}\lambda_i\cdot
\sum_{\begin{subarray}{c}
\mu\;\mbox{\scriptsize s.t.}\;|\mu|=n\;\mbox{\scriptsize and}\;\\
\ell(\mu)=n+1-\ell(\lambda)
\end{subarray}}
\frac{1}{\prod_j m_j(\mu)!}.
\end{eqnarray*}
The sum at extreme right can then be evaluated with the help of the formal power series identity
$$
\sum_{\nu}\frac{x^{\ell(\nu)}y^{|\nu|}}{\prod_i m_i(\nu)!}
=\exp\left(\frac{xy}{1-y}\right)
=\sum_{\ell=0}^\infty\sum_{n=\ell}^\infty \frac{(n-1)!}{(n-\ell)!\ell!(\ell-1)!}
x^\ell y^n,
$$
where the sum on the extreme left is extended over all numerical partitions $\nu$.
The proof is complete. \qed

\subsubsection{Graphical interpretation of the set $\GJdM_n$}\label{subsubsection:GraphGJdM}
Fix $(\theta,\sigma)\in \GJ_n(\theta,\sigma)$ and an element $g\in \dMotz_n(\theta,\sigma)$. 
Let $\Tfrak_n(\theta,\sigma)\in \SV_n$ be a Shabat-Voevodsky tree
from which one recovers the pair $(\theta,\sigma)$.
By \eqref{equation:Interaction2}, the function $g$ 
factors through $\Orbit_n(\sigma)$ and thus may be construed as a function
defined on the set of black vertices of $\Tfrak_n(\theta,\sigma)$.  
In other words, $g$ can be interpreted
as a ``painting over'' of the black vertices of $\Tfrak_n(\theta,\sigma)$ using three new colors, say blue, green, and red, 
corresponding to $0$, $-1$ and $1$, respectively.
Let $\Tfrak_n(\theta,\sigma;g)$ denote the resulting four-colored edge-labeled planar tree of $n$ edges.
Thus we have constructed a bijection 
identifying $\GJdM_n$ with the set of (equivalence classes of) edge-labeled 
vertex-colored planar trees of $n$ edges where the vertex-coloring has to obey the following rules:
\begin{itemize}
\item Only four colors (blue, white, green and red) are used altogether.
\item In any pair of adjacent vertices, exactly one is white.
\item Every white vertex has as many red neighbors as green,
cf. \eqref{equation:Interaction1}.
\item Every green vertex has degree one, i.e., is a leaf,
cf. \eqref{equation:Interaction3}.
\item Every blue vertex has degree two,
cf. \eqref{equation:Interaction3.5} and \eqref{equation:Interaction4}.
\end{itemize}
If we restrict attention  to the set $\GJdM_n(\theta)$,
then the representing trees $\Tfrak$ are required to satisfy the further condition that $\theta=\theta_\Tfrak$.
\begin{Lemma}\label{Lemma:NoBlue}
Fix $(\theta,\sigma,g)\in \GJdM_n$ such that every $\theta$-orbit has even cardinality.
(i) Then $|\{g=-1\}\cap A|=|A|/2$ for every block $A\in \Orbit_n(\theta)$.
(ii) Furthermore,  every element of $\{g=-1\}$ is a fixed point of $\sigma$.
\end{Lemma}
\proof Let $\Tfrak=\Tfrak(\theta,\sigma;g)$.
In graphical language, the claim being made here is that for every white vertex of $\Tfrak$ exactly half of its neighbors are green leaves.
In view of the coloring rules we have only to rule out the existence of blue vertices.
In any case, every white vertex of $\Tfrak$ has an even number of blue neighbors. 
Were $\Tfrak$ to have at least one blue vertex, the coloring rules
would force a circuit to exist, which is a contradiction.
\qed

\subsubsection{The cancellation construction}
Let $X\subset A$ be an inclusion of finite sets. Let $\tau\in S_A$ be a permutation.
For $i\in A\setminus X$, let $\mu(\tau,A,X,i)$ be the least of the positive integers $m$ such that $\tau^m(i)\in A\setminus X$.
We define $\tau\backslash X\in S_{A\setminus X}$ by the formula
$(\tau\backslash X)(i)=\tau^{m(\tau,A,X,i)}(i)$ for $i\in A\setminus X$.
A more intuitively accessible if less precise description of $\tau\backslash X$ is as follows.
Firstly, one writes out the canonical factorization
of $\tau$.
Secondly, one strikes all elements of $X$ from the factorization.
Thirdly and finally, one discards all cycles reduced to length $\leq 1$ by the operation
of striking elements of $X$.
The resulting expression is then the canonical factorization of $\tau\backslash X$.  
\subsubsection{Snipping off green leaves}\label{subsubsection:Snip}
Fix an Eulerian partition $\lambda$ (all parts even)
along with some $\theta\in S_n$ such that $\theta\sim \lambda$.
As usual let $\ell=\ell(\lambda)=\ell(\theta)$.
Fix any set $X\subset \langle n\rangle$ intersecting each block $A\in\Orbit_n(\theta)$
in a set of cardinality $|A|/2$.
Let $\GJdM_n(\theta,X)$ denote the subset of $\GJdM_n(\theta)$
consisting of $(\sigma,g)$ such that $\{g=-1\}=X$.
In order to reconcile
the expression on the right side of
\eqref{equation:TutteBis} with the expression on the right side
of \eqref{equation:MainResult},
it will be enough by Lemma \ref{Lemma:NoBlue}
to prove that
\begin{equation}\label{equation:Nuts}
|\GJdM_n(\theta,X)|=\frac{(\frac{n}{2}-1)!}{(\frac{n}{2}-\ell+1)!}
\prod_{i=1}^\ell \frac{\lambda_i}{2}.
\end{equation}
Now pick $(\sigma,g)\in \GJdM_n(\theta,X)$ arbitrarily
and let $\Tfrak=\Tfrak(\theta,\sigma;g)$.
In turn, let $\Tfrak'$ be the object
obtained from $\Tfrak$ by snipping off each green leaf and attached ``stem,'' while leaving the white vertex at the other end in place, and  blackening all red vertices.
We emphasize that the sets of white vertices of $\Tfrak$ and $\Tfrak'$ are exactly the same. 
Then (the equivalence class of) the object $\Tfrak'$ belongs to $\SV_{\langle n\rangle\setminus X}$
and satisfies $\theta_{\Tfrak'}=\theta\backslash X$.   
Now on the one hand, given $\Tfrak'$, it is clear how to reconstruct $\Tfrak$
by reattaching the green leaves and labeled stems that were snipped off. 
On the other hand, the possible objects $\Tfrak'$  are counted (up to equivalence)
by Lemma \ref{Lemma:GJEnumeration}, and thus  \eqref{equation:Nuts} indeed holds.
In this way the right sides of \eqref{equation:MainResult}
and  \eqref{equation:TutteBis} are  reconciled.
\begin{figure}[hbt]
\centering
\includegraphics[scale=1]{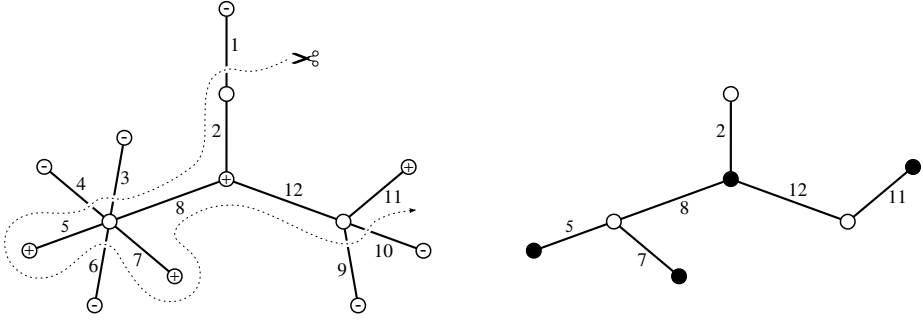}
\caption{
Snipping green leaves off the (left) tree   
$\mathfrak T(\theta,\sigma;g)$ in ${\rm GJdM}_{12}(\theta,X)$ 
where $\theta=(1,2)(3,4,5,6,7,8)(9,10,11,12)$, ${X=\{1,3,4,6,9,10\}}$, 
$\theta\setminus X=(5,7,8)(11,12)$,
and $\sigma=(2,8,12)$, then blackening red vertices, 
gives the (right) Shabat-Voevodsky tree $\Tfrak'$ in ${\rm SV}_{\langle 12\rangle\setminus X}$
for which $\theta_{\Tfrak'}=\theta\setminus X$.
\label{Tree:Surgery:Fig}
}
\end{figure}

\subsection{Comparison of $\GJdM_n$ with the class of labeled mobiles}
\label{subsection:PMM}
In \S\ref{subsubsection:GraphGJdM} we gave a graphical representation for 
members of $\GJdM_n(\theta)$ as certain decorated planar trees $\mathfrak T$.  
In this section we reconcile $\GJdM_n(\theta)$ with some well-known objects   
by giving a bijection between the set of representing trees $\mathfrak T$ 
and a certain subset $\mathscr T_n(\theta)$ of edge-labeled \emph{generalized mobiles} 
considered by Bouttier-Di Francesco-Guitter in~\cite{BouFraGui}.

Fix $(\theta,\sigma)\in \GJ_n(\theta,\sigma)$ and an element $g\in\dMotz_n(\theta,\sigma)$.  
Let $\mathfrak T_n(\theta,\sigma;g)$ be the graphical representation given in \S\ref{subsubsection:GraphGJdM}.  
From the tree we construct a different decorated tree   
$\widehat{\mathfrak T}_n(\theta,\sigma;h)$ through a reversible process.  
Here $h$ denotes the unique function satisfying 
$h\circ \theta - h=g$, $h\circ \sigma=h$ and $\min_{i\in\nn} h(i)=0$; existence and uniqueness 
follow from Proposition~\ref{Proposition:EquiNumerous} 
below.

\begin{figure}[tbh]
\centering
\begin{subfigure}[t]{0.32\textwidth}
\centering\includegraphics{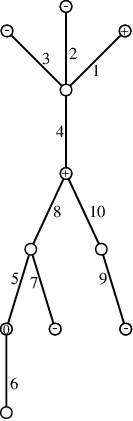}
\caption{}
\label{DJdM:Orig:Fig}
\end{subfigure}
\begin{subfigure}[t]{0.32\textwidth}
\centering\includegraphics{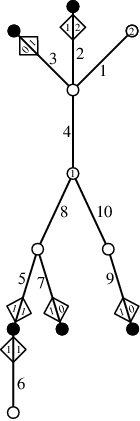}
\caption{}
\label{DJdM:Mobile:Fig}
\end{subfigure}
\caption{
The tree $\mathfrak T_{10}(\theta,\sigma;g)$ 
in (a) 
represents a triple $(\theta,\sigma;g)\in \mathrm{GJdM}_{10}$ 
(from Figure~\ref{Diagrams:Fig} above), 
and the decorated tree in (b) is   
the corresponding labeled mobile $\widehat{\mathfrak T}_{10}(\theta,\sigma;h)$.  
The function $h$ is given by $h(4)=h(8)=h(10)=0$, $h(1)=h(3)=h(5)=h(6)=h(7)=h(9)=1$, and $h(2)=2$.
}\label{DJdM:To:Mobile:Fig}
\end{figure}

To each edge in $\mathfrak T_n(\theta,\sigma;g)$ 
incident to a green or blue vertex,  
attach a ``flag'' pointing in the positive (resp.,~negative) 
direction around the incident white vertex, and write 
the value $h(\theta(i))$ (resp., $h(i)$) on the flag, 
where $i$ is the label of the edge.
Label each red vertex with the value $h(\theta(i))$ for any 
incident edge label $i$.    
Now paint the colored (blue, green, red) vertices black and 
let $\widehat{\mathfrak T}_n(\theta,\sigma;h)$ denote 
the resulting decorated Shabat-Voevodsky tree.  
From $\widehat{\mathfrak T}_n(\theta,\sigma;h)$
we can recover $\mathfrak T_n(\theta,\sigma;g)$
by painting 
labeled vertices red, removing vertex labels, and  
painting each remaining 
black vertex green (resp.,~blue) 
if it has exactly one (resp.,~two) incident edges; see Figure~\ref{DJdM:To:Mobile:Fig}.

Let $\mathscr T_n(\theta)$ denote the set of trees $\widehat{\mathfrak T}_n(\theta,\sigma;h)$ 
constructed from $\GJdM_n(\theta)$.  
Stripping edge labels (but not flags) from the trees gives 
another set $\mathscr T_n(\lambda)$ of trees with the same 
white vertex degree distribution $\lambda\sim \theta$, 
but no edge labels.  Bouttier-Di Francesco-Guitter 
gave an explicit bijection (see~\cite[\S3 and \S4.2]{BouFraGui}) from 
$\mathscr T_n(\lambda)$ onto the set of pairs $(M,v)$  
where $M$ is a planar map with face degree distribution $\lambda$ and $v$ is 
a vertex of $M$, or dually, from $\mathscr T_n(\lambda)$
onto the set of pairs $(M,f)$ where $M$ is a planar map 
with vertex degree distribution $\lambda$ and $f$ is a face of $M$.  
The bijection extends to a labeled version in the natural way
to give a bijection 
from $\mathscr T_n(\theta)$ (and hence from $\GJdM_n(\theta)$) onto the set of all pairs 
$(M,f)$ where $M\in \Map_n(\theta)$ and $f$ is a face of $M$, thus reconciling 
Theorem~\ref{Theorem:MainResult} with~\cite{BouFraGui}.

 Finally we remark that  all the constructions sketched immediately above
 as well as the bijection of \cite{BouFraGui} can be made explicit in a framework emphasizing permutation pairs,
in the spirit of \cite{CoriMachi} and \cite{GouldenJackson}.
This topic will be discussed by the third author on another occasion.

\section{Joint cumulants of functions of a Gaussian random vector}\label{section:BKAR}
Our goal in this section is to derive a delicate expansion of the right side of formula  \eqref{equation:PreMotzkinExpansion}.
(See Proposition \ref{Proposition:ToolCrush} below.) 
We obtain this expansion by specializing a general representation for the joint cumulant of several polynomial functions of a given Gaussian random vector. (See Theorem \ref{Theorem:MainTool} below).
 We obtain the latter representation
by applying the BKAR formula from rigorous statistical mechanics. (See Theorem \ref{Theorem:BKARarboreal} below.)  We have written this section anticipating that the reader would be unfamiliar with the BKAR formula but otherwise familiar with common tools from combinatorics and probability.  Accordingly, we have made our discussion of the BKAR formalism more or less self-contained, if rather compressed. 

 \subsection{Joint cumulants and related apparatus}\label{subsection:JointCumulants}
We briefly recall the formalism of set partitions, M\"{o}bius inversion and joint cumulants, mostly for the purpose of fixing notation.
See \cite{Rota} for the foundations.
(Caution: we do not follow the notation of this reference too closely.)
See also \cite[Section 8.6]{Jacobson} for a short treatment
of generalities concerning M\"{o}bius functions of finite posets.
See  \cite[II.12.8]{Shiryaev} for a probability textbook treatment of joint cumulants.

\subsubsection{Set partitions}
Let $n$ be a positive integer. Recall our abbreviated notation $\langle n\rangle=\{1,\dots,n\}$. 
 A {\em set partition} of $\langle n\rangle$ (or, context permitting, simply a {\em partition}) is by definition a disjoint family of nonempty subsets of $\langle n\rangle$
the union of which equals $\langle n\rangle$. The family of partitions of $\langle n\rangle$ will be denoted by $\Part_n$.
Given $\Pi\in \Part_n$, each member of $\Pi$ is called a {\em block}.
Given $\Pi_1,\Pi_2\in \Part_n$ we write $\Pi_1\leq\Pi_2$ and say that $\Pi_1$ is a {\em refinement}
of $\Pi_2$ if for every block $A\in \Pi_1$ there exists some block $B\in \Pi_2$ such that $A\subset B$.
We also write $\Pi_1<\Pi_2$ if $\Pi_1\leq \Pi_2$ but $\Pi_1\neq \Pi_2$.
Thus partially ordered by refinement, $\Part_n$ becomes a {\em lattice}, i.e., a poset in which  
every family $F$ of elements has a  greatest lower bound $\wedge F$ and a least upper bound $\vee F$.
The least partition $\{\{i\}\mid i\in \langle n\rangle\}=\wedge\Part_n=\vee \emptyset$ will be denoted by $\zero_n$.
The greatest partition $\{\langle n\rangle\}=\vee\Part_n=\wedge \emptyset$  will be denoted by $\one_n$.
For $\Pi_1,\Pi_2\in \Part_n$ such that $\Pi_1\leq \Pi_2$, let 
$$[\Pi_1:\Pi_2]=\{\Pi\in \Part_n\mid \Pi_1\leq \Pi\leq \Pi_2\},$$
which one calls the {\em interval} bounded below by $\Pi_1$ and above by $\Pi_2$.

\subsubsection{The M\"{o}bius function of $\Part_n$}\label{subsubsection:MoebiusReview}
The {\em M\"{o}bius function} 
$$\mu=\mu_{\Part_n}=((\Pi_1,\Pi_2)\mapsto \mu(\Pi_1:\Pi_2)):\Part_n\times \Part_n\rightarrow \ZZ$$
is that function which, when viewed as a $\Part_n$-by-$\Part_n$ matrix,
is inverse to the {\em incidence matrix }
$$((\Pi_1,\Pi_2)\mapsto \indicator{\Pi_1\leq \Pi_2}):\Part_n\times \Part_n\rightarrow \{0,1\}.$$
Since the incidence matrix is upper unitriangular, so also is the matrix $\mu$, i.e., 
$$
\mu(\Pi_1:\Pi_2)=\delta_{\Pi_1,\Pi_2}\;\;\mbox{unless $\Pi_1<\Pi_2$.}
$$
By definition of $\mu$ one has
\begin{equation}\label{equation:MoebiusInversion}
\sum_{\Pi\in [\Pi_1:\Pi_2]}\mu(\Pi:\Pi_2)=\sum_{\Pi\in [\Pi_1:\Pi_2]}\mu(\Pi_1:\Pi)=\delta_{\Pi_1,\Pi_2}\;\;
\mbox{for $\Pi_1\leq \Pi_2$.}
\end{equation}
This is the {\em M\"{o}bius inversion formula} for the lattice $\Part_n$.  
The M\"{o}bius function is given explicitly  for $\Pi_1\leq \Pi_2$ by the expression
\begin{equation}\label{equation:MoebiusFormula}
\mu(\Pi_1:\Pi_2)
=
\prod_{B\in \Pi_2}
(-1)^{|\{A\in \Pi_1\mid A\subset B\}|-1}(|\{A\in \Pi_1\mid A\subset B\}|-1)!.
\end{equation}
(See \cite[Corollary, p. 360]{Rota}.) 
Note that $\mu(\Pi_1:\Pi_2)$ depends only on the isomorphism class of the poset $[\Pi_1:\Pi_2]$.
(This last remark holds for the M\"{o}bius function of any finite poset.)

\subsubsection{The joint cumulant functional} \label{subsubsection:JointCumulantFunctional}
Let $S$ be a finite index set. Let $\{X_i\}_{i\in S}$ be a family of real-valued random variables each member of which has absolute moments of all orders. The {\em joint cumulant} of these variables
is defined by the formula
\begin{equation}\label{equation:BasicCumulantDef}
\kappa\left(\{X_i\}_{i\in S}\right)=
\left(\prod_{i\in S}\frac{\partial}{\partial t_i}\right)
\log \Ebold \exp\left(\sum_{i\in S}t_iX_i\right)\bigg\vert_{\mbox{\scriptsize $t_i=0$ for $i\in S$}},
\end{equation}
where  the variables $t_i$ are treated formally. Hereafter we suppose for simplicity that $\{X_i\}_{i\in S}=\{X_i\}_{i=1}^n$.
Via formula \eqref{equation:MoebiusFormula} one has an equivalent expression
\begin{equation}\label{equation:CumulantDef}
\kappa(X_1,\dots,X_n)=\sum_{\Pi\in \Part_n}
\mu(\Pi:\one_n)\prod_{A\in \Pi}\Ebold \prod_{i\in A}X_i
\end{equation} 
for the joint cumulant functional.
By the M\"{o}bius inversion formula \eqref{equation:MoebiusInversion} one then has
an expansion
\begin{equation}\label{equation:CumulantDefBis}
\Ebold \prod_{i=1}^nX_i=\sum_{\Pi\in \Part_n}\;
\prod_{A\in \Pi}\kappa\left(\{X_i\}_{i\in A}\right).
\end{equation}

\subsubsection{The Isserlis-Wick formula}
Suppose now $X_1,\dots,X_n$ are real random variables with a centered Gaussian joint distribution.
Using \eqref{equation:BasicCumulantDef} one can show   that
the joint cumulant of three or more random variables with a Gaussian joint distribution vanishes identically.
Thus, after substituting into \eqref{equation:CumulantDefBis}, one obtains 
the relation
\begin{equation}\label{equation:WickFormula}
\Ebold \prod_{i=1}^nX_i=\sum_{\begin{subarray}{c}
\Pi\in \Part_n\;\mbox{\scriptsize s.t.}\\
\mbox{\scriptsize all blocks are}\\
\mbox{\scriptsize of cardinality $2$}
\end{subarray}}\prod_{\{i,j\}\in \Pi}\Ebold X_iX_j,
\end{equation}
 known as the {\em Wick formula} among physicists but in fact due to Isserlis
\cite{Isserlis}.

\subsubsection{Trivial generalization of \eqref{equation:CumulantDef}}
Suppose that for some partition $\Theta\in \Part_n$ one is given a family of real random variables
$\{Y_A\}_{A\in \Theta}$ with absolute moments of all orders.
Then by \eqref{equation:MoebiusFormula} one has for  formula \eqref{equation:CumulantDef} a trivial
generalization
\begin{equation}\label{equation:CumulantDefTer}
\kappa\left(\{Y_A\}_{A\in \Theta}\right)=\sum_{\Pi\in [\Theta:\one_n]}\;
\mu(\Pi:\one_n)\prod_{B\in \Pi}\Ebold \prod_{\begin{subarray}{c}
A\in \Theta\\
\mbox{\scriptsize s.t.}\;A\subset B
\end{subarray}}Y_A
\end{equation}
which will be especially important in the sequel.

\subsection{The probability measures
$\PP_\Gamma^\Theta$}
Under this heading we make the key definition figuring in the BKAR formula and
(hence) in the statement of the refined expansion of the right side of \eqref{equation:PreMotzkinExpansion} we are aiming to obtain.
We actually give a couple of equivalent definitions, each of which has its uses.
The formalism we set up here will be in use throughout the paper.

\subsubsection{The set $\Bond_n$}
Let 
$$\Bond_n=\{\{i,j\}\subset \langle n\rangle\mid i,j\in \langle n\rangle,\;i\neq j\}\subset 2^{\langle n\rangle}.$$
For each subset $\Gamma\subset \Bond_n$ and partition $\Theta\in \Part_n$, with some abuse of notation, let $\Gamma\vee\Theta$ denote the greatest lower bound of the family of partitions $\Psi\in [\Theta:\one_n]$
such that each member of $\Gamma$ is contained in some block of $\Psi$. 
Roughly speaking $\Gamma\vee\Theta$ arises from $\Theta$ by coalescing pairs of blocks whenever they are ``bonded'' by some member of $\Gamma$. 

\begin{figure}[hbt]
\centering
\includegraphics[scale=1]{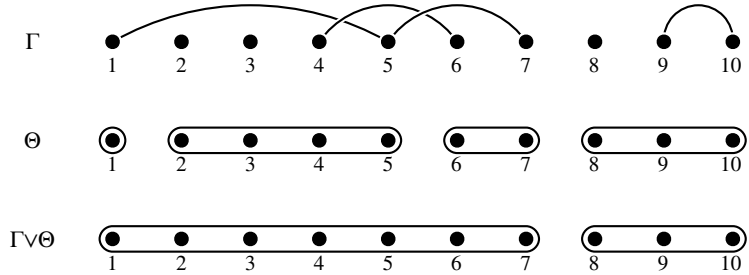}
\caption{
The set of bonds $\Gamma=\{\{1,5\},\{4,6\},\{5,7\},\{9,10\}\}\subset \Bond_{10}$, 
the partition  
$\Theta=\{\{1\},\{2,3,4,5\},\{6,7\},\{8,9,10\}\}\in \Part_{10}$, 
and the partition $\Gamma\vee \Theta=\{\{1,2,3,4,5,6,7\},\{8,9,10\}\}\in \Part_{10}$.
\label{Join:Fig}
}
\end{figure}

\subsubsection{The graphs $\Gfrak(\Theta,\Gamma)$}
Given $\Theta\in \Part_n$  and $\Gamma\subset \Bond_n$,
we define a graph $\Gfrak(\Theta,\Gamma)$ by the following conventions:
\begin{itemize}
\item Each member of $\Theta$ is interpreted as a vertex.
\item Each member of $\Gamma$ is interpreted as an edge.
\item For all edges $e=\{i,j\}\in \Gamma$ and vertices $A,B\in \Theta$ such that $i\in A$ and $j\in B$,
the set of endpoints of $e$ is declared to be $\{A,B\}$.
\end{itemize}
The graph $\Gfrak(\Theta,\Gamma)$ has
in general multiple edges and loops joining a vertex to itself. 
Most graphs we need to consider in this paper arise naturally in the form $\Gfrak(\Theta,\Gamma)$.
Note that the family of connected components of the graph $\Gfrak(\Theta,\Gamma)$
is canonically in bijection with the set $\Gamma\vee\Theta$.

\begin{figure}[hbt]
\centering
\includegraphics[scale=1]{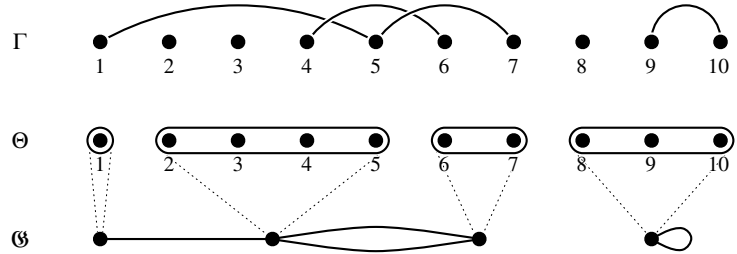}
\caption{
The set of bonds $\Gamma=\{\{1,5\},\{4,6\},\{5,7\},\{9,10\}\}\subset \Bond_{10}$, 
the set partition  
$\Theta=\{\{1\},\{2,3,4,5\},\{6,7\},\{8,9,10\}\}\in \Part_{10}$, 
and the graph $\mathfrak G=\mathfrak G(\Theta,\Gamma)$.
\label{BondGraphDot:Fig}
}
\end{figure}

\subsubsection{The set $\Tree_n(\Theta)$}
For $\Theta\in \Part_n$, let $\Tree_n(\Theta)$ denote the set of $\Gamma\subset \Bond_n$
such that $\one_n=\Gamma\vee\Theta$
and $|\Gamma|+1=|\Theta|$.
Equivalently, $\Tree_n(\Theta)$ is the set whose members are sets $\Gamma\subset \Bond_n$
such that the graph $\Gfrak(\Theta,\Gamma)$ with vertex set $\Theta$ and edge set $\Gamma$ is connected
and has Euler characteristic $1$, i.e., is a tree. 

\begin{figure}[H]
\centering
\includegraphics[scale=1]{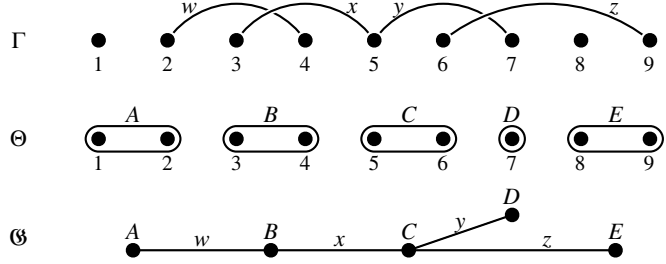}
\caption{
A member $\Gamma=\{\{2,4\},\{3,5\},\{5,7\},\{6,9\}\}\subset \Bond_9$ 
of ${\rm Tree}_9(\Theta)$ for $\Theta=\{\{1,2\},\{3,4\},\{5,6\},\{7\},\{8,9\}\}\in \Part_9$,  
and the 
 tree $\mathfrak G=\mathfrak G(\Theta,\Gamma)$.
\label{Tree:Fig}
}
\end{figure}
\subsubsection{The Schur Product Theorem and related notation}\label{subsection:MatrixNotation}
\begin{itemize}
\item Let $\Mat_n$ denote the space of $n$-by-$n$ matrices with real entries
\item Let $\Sym_n\subset \Mat_n$ denote the space of symmetric matrices.
\item Let $\Sym_n^+=\{Q\in \Sym_n\mid \mbox{$Q$ is positive semidefinite}\}$.
\item For $A,B\in \Mat_n$, recall that the {\em Hadamard product} 
(alternatively and arguably more correctly: the {\em Schur product}) $A\star B\in \Mat_n$ is defined by the formula $(A\star B)(i,j)=A(i,j)B(i,j)$
 (entry-by-entry multiplication).
\end{itemize}
According to the {\em Schur Product Theorem} \cite{Schur}
 if $A,B\in \Sym_n^+$, then
$A\star B\in \Sym_n^+$. The latter fact is of extreme importance in the sequel.
\subsubsection{The set $\Qfrak_n$}
Let $\Qfrak_n$ denote the set consisting of all matrices $Q\in \Sym_n^+$
with the following properties:
\begin{itemize}
\item All entries of $Q$ belong to the closed unit interval $[0,1]$.
\item All diagonal entries of $Q$ are equal to $1$.
\end{itemize}
Since the $n$-by-$n$ identity matrix belongs to $\Qfrak_n$,  the latter set is not empty.
It is easy to see that the set $\Qfrak_n$ is closed, convex, bounded and hence compact.
 For $\Theta\in \Part_n$ and $\Gamma\in \Tree_n(\Theta)$,
the probability measure $\PP^\Theta_\Gamma$ we aim to define
will be defined on the set $\Qfrak_n$.

\subsubsection{The matrix representation of partitions}
Given $\Pi\in \Part_n$, we define the matrix $[\Pi]\in\Mat_n$ to have entries
$$[\Pi](i,j)=\left\{\begin{array}{rl}1&\mbox{if $i$ and $j$ belong to the same block of $\Pi$,}\\
0&\mbox{otherwise.}
\end{array}\right.
$$
We say that the matrix $[\Pi]$ thus defined {\em represents} $\Pi$. 
For example, one has
$$
[\{\{1,2,3\},\{4,5\}\}]=
\left[\begin{array}{ccccccccc}
1&1&1&0&0\\
1&1&1&0&0\\
1&1&1&0&0\\
0&0&0&1&1\\
0&0&0&1&1
\end{array}\right].
$$
In particular, $[\zero_n]$ is the $n$-by-$n$ identity matrix
and $[\one_n]$ is the $n$-by-$n$ matrix with all entries equal to $1$.
Finally and crucially, note that $[\Pi]\in \Qfrak_n$ for $\Pi\in \Part_n$.

\subsubsection{Definition of $\PP_\Gamma^\Theta$} 
Fix $\Theta\in \Part_n$ arbitrarily and let $k=|\Theta|$.
We define a family 
$$\left\{\PP_\Gamma^\Theta\right\}_{\Gamma\in \Tree_n(\Theta)}$$ of probability measures
on $\Qfrak_n$ by requiring  the integration formula
\begin{eqnarray}\label{equation:UndergraduatePP}
&&\sum_{\Gamma\in \Tree_n(\Theta)}\int f_\Gamma \,\dd \PP_\Gamma^\Theta\\
\nonumber&=&\sum_{\begin{subarray}{c}
(e_1,\dots,e_{k-1})\in \Bond_n^{k-1}\\
\mbox{\scriptsize s.t.}\,\{e_1,\dots,e_{k-1}\}\in \Tree_n(\Theta)
\end{subarray}}\begin{array}{c}
\\
\displaystyle\int\cdots \int\\
\scriptstyle 1=t_0>t_1>\cdots>t_{k-1}>t_k=0
\end{array}\\
\nonumber&&f_{\{e_1,\dots,e_{k-1}\}}\left(\sum_{\alpha=0}^{k-1}(t_{\alpha}-t_{\alpha+1})[\{e_{1},\dots,e_{\alpha}\}\vee \Theta]
\right)\prod_{\alpha=1}^{k-1}\dd t_\alpha\\
\nonumber&=&\sum_{\begin{subarray}{c}
(e_1,\dots,e_{k-1})\in \Bond_n^{k-1}\\
\mbox{\scriptsize s.t.}\,\{e_1,\dots,e_{k-1}\}\in \Tree_n(\Theta)
\end{subarray}}
\int_0^1\dd t_1\int_0^{t_1}\dd t_2\cdots \int_0^{t_{k-2}}\dd t_{k-1}\\
\nonumber&&f_{\{e_1,\dots,e_{k-1}\}}\left([\Theta]+\sum_{\alpha=1}^{k-1}t_\alpha
([\{e_1,\dots,e_\alpha\}\vee \Theta]-
[\{e_1,\dots,e_{\alpha-1}\}\vee \Theta])\right)
\end{eqnarray}
to hold for every family
 $$\{f_\Gamma:\Qfrak_n\rightarrow\RR\}_{\Gamma\in \Tree_n(\Theta)}$$
 of continuous functions.

\begin{Lemma}[Alternate characterization of $\PP^\Theta_\Gamma$]\label{Lemma:PTreeChar}
Fix $\Theta\in \Part_n$ and $\Gamma\in \Tree_n(\Theta)$. Let $k=|\Theta|$.
Let $X\in \Qfrak_n$ be a random matrix with law $\PP_\Gamma^\Theta$.
For $\{i,j\}\in \Bond_n$ and blocks $A,B\in \Theta$ such that $i\in A$ and $j\in B$,
let $\Gamma(i,j)\subset \Gamma$ be the subset consisting of edges visited by the
unique geodesic walk in the tree $\Gfrak(\Theta,\Gamma)$ joining $A$ to $B$.  Then the following
statements concerning the random matrix $X$ hold:\\
\begin{enumerate}
\item[(i)]
The family of matrix entries 
$$\{X(i,j)\mid 1\leq i<j\leq n\;\mbox{and}\; \{i,j\}\in \Gamma\}$$ is i.i.d.uniformly distributed in $(0,1)$.\\
\item[(ii)]
For all $i,j\in \langle n\rangle$ one has 
$$X(i,j)=\min\left(\{1\}\cup\{X(i',j')\mid \{i',j'\}\in \Gamma(i,j)\}\right)$$ almost surely.
\end{enumerate}
\end{Lemma}
\noindent The lemma reconciles the definition \eqref{equation:UndergraduatePP} of $\PP^\Theta_\Gamma$ given 
above with the form of the definition typical in the literature. 
Later this lemma will permit
us to calculate certain integrals coming up in the proof of Theorem \ref{Theorem:MainResult}.
(See Proposition \ref{Proposition:LastCut} below.)
As a point of contact with the literature, we mention \cite[Lemma 4.1]{Faris} which in the context of the matroidal generalization of the BKAR formula serves an end very similar to that served by Lemma \ref{Lemma:PTreeChar}.

\proof
We begin by building an explicit random matrix with law $\PP_\Gamma^\Theta$.
Let
$$T=(T_1,\dots,T_{k-1})$$ be a real random vector uniformly distributed in the simplex
$$\{(t_1,\dots,t_{k-1})\in \RR^{k-1}\mid 1>t_1>\cdots>t_{k-1}>0\}.$$
For convenience let $T_0=1$ and $T_k=0$. 
Write $\Gamma=\{e_1,\dots,e_{k-1}\}$.
Let $\rho\in S_{k-1}$ 
be a uniformly distributed random permutation independent of $T$.
Consider the random matrix
\begin{eqnarray}\label{equation:Xformula}
&&[\Theta]+\sum_{\alpha=1}^{k-1}T_\alpha([\{e_{\rho(1)},\dots,e_{\rho(\alpha)}\}\vee \Theta]-
[\{e_{\rho(1)},\dots,e_{\rho(\alpha-1)}\}\vee \Theta])\\
\nonumber&=&\sum_{\alpha=0}^{k-1}(T_{\alpha}-T_{\alpha+1})[\{e_{\rho(1)},\dots,e_{\rho(\alpha)}\}\vee \Theta]
\end{eqnarray}
which clearly takes its values in $\Qfrak_n$. It is a trivial matter to confirm that
the law on $\Qfrak_n$ of the random matrix \eqref{equation:Xformula}
is $\PP_\Gamma^\Theta$.  Without loss of generality we may identify the given random matrix $X$
with the random matrix  \eqref{equation:Xformula}.

Fix $\{i,j\}\in \Bond_n$ and blocks $A,B\in \Theta$ such that $i\in A$ and $j\in B$.
 It will be enough to evaluate the matrix entry $X(i,j)$ in terms of $T$ and $\rho$. 
We begin by observing that $X(i,j)=T_\beta$ where the (random) index $\beta$ is the least index $\alpha$ such that
$$[\{e_{\rho(1)},\dots,e_{\rho(\alpha)}\}\vee \Theta](i,j)=1.$$ 
Equivalently, $\beta$ is the least index $\alpha$ such that $A$ and $B$ are connected by some walk in the (random) forest
$$\Gfrak(\Theta,\{e_{\rho(1)},\dots,e_{\rho(\alpha)}\}).$$
Now $e\in \Gamma$ satisfies $e\in \Gamma(i,j)$ if and only if $A$
and $B$ are NOT joined by a walk in the forest $\Gfrak(\Theta,\Gamma\setminus \{e\})$.
Thus $\beta$ is the least index $\alpha$ such that
$$\{e_{\rho(1)},\dots,e_{\rho(\alpha)}\}\supset \Gamma(i,j).$$
By this reasoning we arrive at the formula
$$X(i,j)=\min(\{1\}\cup\{T_{\rho^{-1}(\alpha)}\mid \alpha=1,\dots,k-1\;\mbox{s.t.}\;e_\alpha\in \Gamma(i,j)\}).$$
Now write $e_i=\{a_i,b_i\}$ where $a_i<b_i$ for $i=1,\dots,k-1$. It is clear that the random vector 
$$(X(a_1,b_1),\dots,X(a_{k-1},b_{k-1}))=(T_{\rho^{-1}(1)},\dots,T_{\rho^{-1}(k-1)})$$
is uniformly distributed in the cube $(0,1)^{k-1}$.  Statements (i) and (ii) follow.
\qed

\begin{figure}[hbt]
\centering
\includegraphics[scale=1]{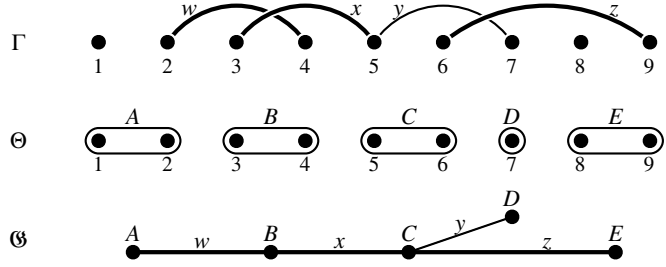}
\caption{
For $\Gamma=\{w,x,y,z\}\subset \Bond_9$ and 
$\Theta=\{A,B,C,D,E\}\in \Part_9$ as drawn above, 
$\Gamma(1,8)=\Gamma(2,8)=\Gamma(1,9)=\Gamma(2,9)=\{w,x,z\}$ corresponds to the geodesic path 
joining $A$ with $E$ in the tree $\mathfrak G=\mathfrak G(\Theta,\Gamma)$.
}
\end{figure}

\subsection{A variant of the BKAR formula}
\subsubsection{Differentiation of functions on $\Sym_n$}
For short, we say that a function \linebreak $f:\Sym_n\rightarrow\RR$
is {\em polynomial} if $f(Q)$ is a polynomial with real coefficients
in the entries of $Q$.
For $e=\{i,j\}\in \Bond_n$, $Q\in \Sym_n$
and polynomial functions $f:\Sym_n\rightarrow\RR$,  let 
$$(\partial_ef)(Q)=\frac{d}{dt}f(Q+t(e_{ij}+e_{ji}))\bigg\vert_{t=0}\;\;(\{e_{\alpha\beta}\}_{\alpha,\beta=1}^n:\;\mbox{standard basis of $\Mat_n$}),$$
thus defining a first order linear differential operator $\partial_e$
acting on polynomial functions defined on $\Sym_n$. 
More generally, for each $\Gamma\subset\Bond_n$ let
$$\partial^\Gamma=\prod_{e\in \Gamma}\partial_e.$$
We have then the following fundamental integration identity.
\begin{Theorem}[Variant of the BKAR  formula]\label{Theorem:BKARarboreal}
For set partitions $\Theta\in \Part_n$
and polynomial functions $f:\Sym_n\rightarrow\RR$ 
we have
\begin{equation}\label{equation:ConnectedBKAR}
\sum_{\Pi\in [\Theta:\one_n]}
\mu(\Pi:\one_n)f([\Pi])=\sum_{\Gamma\in \Tree_n(\Theta)}\int\partial^\Gamma f\; \dd \PP^\Theta_\Gamma.
\end{equation}
\end{Theorem}
\noindent  Formula \eqref{equation:ConnectedBKAR} is true for more general functions $f$ than polynomials ones, but here, for simplicity, we stick to the polynomial case. No greater generality will be needed.
In any case, extension of \eqref{equation:ConnectedBKAR} to larger classes of functions can easily enough be accomplished by polynomial approximation.  For the reader's convenience we supply a short proof of \eqref{equation:ConnectedBKAR} in \S\ref{subsubsection:BKARproof} below; the effort of the setup above renders the proof more or less trivial.

\subsubsection{Background and references concerning the BKAR formula}
\label{subsubsection:BKARbackground}
We mention first of all the paper  \cite{BrydgesK}
of Brydges and Kennedy.
Next we mention the papers  \cite{AbdesselamR1} and \cite{AbdesselamR2}
of Abdesselam and Rivasseau. This explains the abbreviation BKAR.
The notes  \cite{AbdesselamNote} give an accessible introduction
to the BKAR formula and many further references.
  The paper \cite{AbdesselamPS} is a typical application of the BKAR formula wherein the latter is used to bound joint cumulants.  The recent paper \cite{Faris} generalizes the BKAR formula in a natural way
to the setting of matroids. 
The BKAR formula is a relatively recent development in a very old and well-established line of research in statistical mechanics focused on {\em cluster expansions}. Concerning the vast literature of cluster expansions, 
we refer the reader to \cite{BrydgesH}, \cite{BrydgesM}
and \cite{FarisCluster} as possible entry points.

\subsubsection{Example} 
\begin{figure}[hbt]
\centering
\includegraphics[scale=0.8]{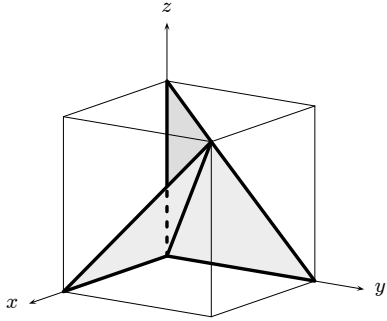}
\caption{This is a depiction of the set of points $(x,y,z)\in [0,1]^3$ such that
$x=y\leq z$ or $y=z\leq x$ or $z=x\leq y$.
\label{FinSt.fig}}
\end{figure}

The case $\Theta=\zero_3\in \Part_3$ of \eqref{equation:ConnectedBKAR} boils down to the
following sophomore calculus exercise:
\begin{eqnarray}\label{equation:Part3FTC}
\\
\nonumber&&\int_0^1\dd u \int_0^u \dd v \,[(f_{xy}+f_{xz})(u,v,v)
+(f_{yx}+f_{yz})(v,u,v)+(f_{zx}+f_{zy})(v,v,u))]\\
\nonumber&=&f(1,1,1)-f(1,0,0)-f(0,1,0)-f(0,0,1)+2f(0,0,0).
\end{eqnarray}
For the purpose of comparison we note the joint cumulant formula
$$\kappa(X,Y,Z)=\Ebold XYZ-\Ebold X \cdot\Ebold YZ-\Ebold Y\cdot \Ebold XZ-\Ebold Z\cdot \Ebold XY+2\Ebold X
\cdot \Ebold Y\cdot\Ebold Z.$$
Figure \ref{FinSt.fig} depicts the set on which all integrations on the left side of formula \eqref{equation:Part3FTC}
are taking place. This set admits interpretation as a geometric realization of the simplicial complex the simplices of which are the chains in the poset $\Part_3$, with the origin corresponding to $\zero_3$
and the point $(1,1,1)$ corresponding to $\one_3$.

\subsubsection{A proof of Theorem \ref{Theorem:BKARarboreal}}\label{subsubsection:BKARproof}
 We may assume that $f$ takes the form
$$
f(Q)=\prod_{1\leq i\leq j\leq n}Q(i,j)^{\nu(i,j)}\;\;\mbox{for $Q\in \Sym_n$}
$$
where
$$\nu=\{\nu(i,j)\mid1\leq i\leq j\leq n\}$$
is a family of nonnegative integers.
Let
$$\supp\,\nu=\{\{i,j\}\in \Bond_n\mid 1\leq i<j\leq n\;\mbox{s.t.}\;\nu(i,j)>0\}.$$
We have
$$\Pi\in [\Theta:\one_n]\Rightarrow f([\Pi])=\indicator{(\supp\,\nu)\vee\Theta\leq \Pi}
$$
and hence by the M\"{o}bius inversion formula \eqref{equation:MoebiusInversion}
we have
\begin{equation}\label{equation:ClinchFTC1}
(\mbox{LHS of \eqref{equation:ConnectedBKAR}})
=\indicator{(\supp\,\nu)\vee\Theta=\one_n}=
\left\{\begin{array}{rl}
1&\mbox{if $\Gfrak(\Theta,\supp\,\nu)$ is connected,}\\
0&\mbox{otherwise.}
\end{array}\right.
\end{equation}
For $\Pi\in \Part_n$ let
$$
N(\Pi)
=\sum_{\begin{subarray}{c}
1\leq i<j\leq n\\
\mbox{\scriptsize $\{i,j\}$ is contained }\\
\mbox{\scriptsize in no block of $\Pi$}
\end{subarray}}\nu(i,j).
$$
Note that for $1\leq i'<j'\leq n$ we have
$$
(\partial_{\{i',j'\}}f)(Q)=
\left\{\begin{array}{cl}
\displaystyle\nu(i',j')\prod_{1\leq i<j\leq n}
Q(i,j)^{\nu(i,j)-\delta_{ii'}\delta_{jj'}}&\mbox{if $\{i',j'\}\in \supp\,\nu$,}\\
0&\mbox{if $\{i',j'\}\not\in \supp\,\nu$.}
\end{array}\right.
$$
Substituting directly into the definition \eqref{equation:UndergraduatePP} of $\PP_\Gamma^\Theta$ we then have
\begin{eqnarray}\label{equation:ClinchFTC2}
&&(\mbox{RHS of \eqref{equation:ConnectedBKAR}})\\
\nonumber&=&\sum_{\begin{subarray}{c}
(e_1,\dots,e_{k-1})\in \Bond_n^{k-1}\\
\mbox{\scriptsize s.t.}\,\{e_1,\dots,e_{k-1}\}\in \Tree_n(\Theta)
\end{subarray}}
\int_0^1\dd t_1\int_0^{t_1}\dd t_2\cdots \int_0^{t_{k-2}}\dd t_{k-1}\\
\nonumber&&\partial_{e_1}\cdots \partial_{e_{k-1}}f\left([\Theta]+\sum_{\alpha=1}^{k-1}t_\alpha(
[\{e_1,\dots,e_\alpha\}\vee \Theta]
-[\{e_1,\dots,e_{\alpha-1}\}\vee\Theta]
)\right)\\
\nonumber&=&\sum_{\begin{subarray}{c}
e=(e_1,\dots,e_{k-1})\in (\supp\,\nu)^{k-1}\\
\mbox{\scriptsize s.t.}\,\{e_1,\dots,e_{k-1}\}\in \Tree_n(\Theta)
\end{subarray}}\int_0^1\dd t_1\int_0^{t_1}\dd t_2\cdots \int_0^{t_{k-2}}\dd t_{k-1}\\
\nonumber&&
\;\;\;\;\;\;\;\;\;\;\;\;\;\;\;\;\;\;\;\;\;\;\;\;\;\; \prod_{\alpha=1}^{k-1}\nu(e_\alpha)t_\alpha^{N(\{e_1,\dots,e_{\alpha-1}\}\vee \Theta)-N(\{e_1,\dots,e_\alpha\}\vee \Theta)-1}\\
\nonumber&=&
\sum_{\begin{subarray}{c}
e=(e_1,\dots,e_{k-1})\in (\supp\, \nu)^{k-1}\\
\mbox{\scriptsize s.t.}\,\{e_1,\dots,e_{k-1}\}\in \Tree_n(\Theta).
\end{subarray}}\;\;\prod_{\alpha=1}^{k-1} \frac{\nu(e_\alpha)}{N(\{e_1,\dots,e_{\alpha-1}\}\vee\Theta)}.
\end{eqnarray}
Now  the right sides of
\eqref{equation:ClinchFTC1}
and \eqref{equation:ClinchFTC2} both vanish if $(\supp\,\nu)\vee \Theta\neq \one_n$.
Otherwise, the right side of \eqref{equation:ClinchFTC1} equals $1$
and thus equals the right side of \eqref{equation:ClinchFTC2} by the lemma recalled immediately below. \qed

\begin{Lemma}
As above, let $\Theta\in \Part_n$ be a set partition and let $k=|\Theta|$.
For every sequence $(e_1,\dots,e_{k-1})\in \Bond_n^{k-1}$
one has $\{e_1,\dots,e_{k-1}\}\in \Tree_n(\Theta)$ if and only if for
$\alpha=1,\dots,k-1$ the set $e_\alpha$ is contained in no block of the set partition $\{e_1,\dots,e_{\alpha-1}\}\vee \Theta$.
\end{Lemma}
\noindent We can safely omit the proof.

\subsection{Formulation of the main technical result}
\subsubsection{Variables}
Let $n$ and $N$ be positive integers.
Let  $$z=\{\{z_{ij}\}_{i=1}^n\}_{j=0}^{2N}$$
be a family of  independent (commutative) algebraic variables. 
Let $\RR[z]$ be the polynomial algebra these variables generate over the real numbers.
We remark that the index $j$ runs here from $0$ to $2N$ rather than, say, from $1$ to $N$ in order
to accommodate the intended application with no adjustment of notation.

\subsubsection{Differential operators}
For $e=\{i,i'\}\in \Bond_n$ we define a partial differential operator
\begin{equation}
D_e=\sum_{j=0}^{2N}\frac{\partial^2}{\partial z_{ij}\partial z_{i'j}}
\end{equation}
acting on the polynomial algebra $\RR[z]$. 
Given $\Gamma\subset \Bond_n$,
we in turn define
\begin{equation}
D^\Gamma=\prod_{e\in \Gamma}D_e.
\end{equation}

\subsubsection{A family of polynomials}
Fix $\Theta\in \Part_n$.
 For each $A\in \Theta$ fix a polynomial

$$f_A\in \RR[\{\{z_{ij}\}_{i\in A}\}_{j=0}^{2N}]\subset \RR[z]$$
and let 
$$f=\prod_{A\in \Theta} f_A\in \RR[z].$$

\subsubsection{Gaussian random variables}Let
$$\zeta=\{\{\zeta_{ij}\}_{i=1}^n\}_{j=0}^{2N}$$
be a family of real random variables with a centered Gaussian joint distribution
such that
\begin{equation}
\Ebold \zeta_{ij}\zeta_{i'j'}=\delta_{jj'}\Ebold \zeta_{i0}\zeta_{i'0}.
\end{equation}
Note that the random vector $\zeta$
has the structure of a family of $2N+1$ i.i.d. copies of the random vector $\{\zeta_{i0}\}_{i=1}^n$. 
But also note that we do not place any restrictions
on the covariances $\Ebold \zeta_{i0}\zeta_{i'0}$. The latter freedom is crucial for the intended application.

\subsubsection{The $Q$-recoupling construction}
For each $Q\in \Qfrak_n$ let
$$\zeta\star Q=\left\{\left\{(\zeta\star Q)_{ij}\right\}_{i=1}^n\right\}_{j=0}^{ 2N}$$
be a family of real random variables with a centered Gaussian joint distribution
characterized by the covariances
\begin{equation}
\Ebold (\zeta\star Q)_{ij}(\zeta\star Q)_{i'j'}=Q(i,i')\Ebold \zeta_{ij}\zeta_{i'j'}=\delta_{jj'}Q(i,i')\Ebold \zeta_{i0}\zeta_{i'0}.
\end{equation}
Such a family $\zeta\star Q$ exists and has a uniquely determined law because the requisite positive-semidefiniteness is guaranteed  by the Schur Product Theorem reviewed in \S\ref{subsection:MatrixNotation} above.  We say that $\zeta\star Q$ arises  from $\zeta$
by {\em $Q$-recoupling}. The probability space on which $\zeta\star Q$ is defined is allowed to depend on $Q$;
it is of no concern to us.
Note that $\zeta\star Q$ has the structure of $2N+1$ i.i.d. of copies of the random vector 
$\{(\zeta\star Q)_{i0}\}_{i=1}^n$.
Note also that for $i=1,\dots,n$ the subfamily $\left\{(\zeta\star Q)_{ij}\right\}_{j=0}^{2N}$ of $\zeta\star Q$ consists
of $2N+1$ i.i.d copies the random variable $\zeta_{i0}$.
Finally, note that $\zeta\star [\one_n]$ is a copy of $\zeta$.

Here is the main technical result of the paper.

\begin{Theorem}\label{Theorem:MainTool}
Notation and assumptions are as above.
We have
\begin{equation}\label{equation:MainTool}
\kappa\left(\left\{f_A(\zeta)\right\}_{A\in \Theta}\right)=
\sum_{\Gamma\in \Tree_n(\Theta)}\;\;
\left(\prod_{\{i,i'\}\in \Gamma}\Ebold \zeta_{i0}\zeta_{i'0}\right)\int\Ebold\left[(D^\Gamma f)(\zeta\star {Q})\right]\;\PP^\Theta_\Gamma(\dd Q).
\end{equation}
\end{Theorem}
\noindent  The proof of \eqref{equation:MainTool} will be given in \S\ref{subsection:ProofOfMainTool} below.
Now to make sense of the right side of \eqref{equation:MainTool} it is necessary to give a consistent
interpretation to expressions of the form
\begin{equation}\label{equation:OuterIntegral}
\int\Ebold\left[g(\zeta\star {Q})\right]\;\PP^\Theta_\Gamma(\dd Q)\;\;\;(g\in \RR[z],\Theta\in \Part_n,\Gamma\in \Tree_n(\Theta)).
\end{equation}
Our convention is invariably to interpret the expressions of form \eqref{equation:OuterIntegral} as iterated integrals.
This interpretation makes sense and indeed yields a well-defined numerical value
because the inner integral $\Ebold\left[g(\zeta\star {Q})\right]$ 
by the Isserlis-Wick formula \eqref{equation:WickFormula} depends polynomially on $Q$.

The next lemma amplifies the theorem by pointing out 
cases in which terms on the right side of \eqref{equation:MainTool} are forced to vanish.
The lemma is the chief means by which we will get the benefit of the theorem in the application.

\begin{Lemma}[``Culling rules'']\label{Lemma:CullingRules} We continue in the setup of Theorem \ref{Theorem:MainTool}.
Fix $\Gamma\in \Tree_n(\Theta)$. Let 
$$Z=\prod_{i=1}^n\prod_{j=0}^{2N}z_{ij}^{\nu_{ij}}\in \CC[z]$$ be a monomial.
For $i\in \langle n\rangle$ let  
$$
\nu_i=\sum_{j=0}^{2N} \nu_{ij}\;\;\mbox{and}\;\;
d_i=|\{e\in \Gamma\mid i\in e\}|.
$$
Then the following statements hold:
\begin{eqnarray}\label{equation:FirstCull}
\prod_{\{i,i'\}\in \Gamma}\sum_{j=0}^{2N}\nu_{ij}\nu_{i'j}=0&\Rightarrow&D^\Gamma Z=0.\\
\label{equation:SecondCull}
\max_{i\in \langle n\rangle}(d_i-\nu_i)>0&\Rightarrow& D^\Gamma Z=0.
\end{eqnarray}
\end{Lemma}
\noindent 
\proof 
Under the hypothesis of \eqref{equation:FirstCull}
there exists some $e\in\Gamma$ such that $D_e Z=0$
and {\em a fortiori} $D^\Gamma Z=0$. Thus \eqref{equation:FirstCull} holds.
Let 
$$\vec{\Gamma}=\{(i,i')\in \langle n\rangle^2\mid \{i,i'\}\in \Gamma\},$$
noting that $|\vec{\Gamma}|=2|\Gamma|=2|\Theta|-2$.
We have a general expansion
\begin{equation}
D^\Gamma
=\sum_{
\mathbf{j}:\Gamma\rightarrow \{0,\dots,2N\}
}\;\;\prod_{(i,i')\in \vec{\Gamma}}\;\frac{\partial}{\partial z_{i,\mathbf{j}(\{i,i'\})}}
\end{equation}
which proves \eqref{equation:SecondCull}.
\qed

\subsubsection{Notes and references}
Formula \eqref{equation:MainTool} hypergeneralizes explicit identities used to prove the Poincar\'{e} inequality
for Gaussian random variables.
For discussion of identities of the latter type see \cite{BobkovGotzeHoudre}, where similar identities
for Bernoulli random variables are  also discussed. 
An analogue of \eqref{equation:MainTool} for Bernoulli random variables would be of considerable interest.

\subsection{Proof of Theorem \ref{Theorem:MainTool}}\label{subsection:ProofOfMainTool}
The next lemma rewrites the Isserlis-Wick formula in a more convenient form involving differential operators.
\begin{Lemma}
\label{Lemma:FancyWickFormula} 
Let $S$ be a finite index set.
Let $t=\{t_i\}_{i\in S}$ be a family of independent commuting algebraic variables.
Each variable $t_i$ is assigned the degree $1$, and we consider the polynomial ring $\RR[t]$
graded by degree.
Let $\tau=\{\tau_i\}_{i\in S}$ be a family of real random variables with a centered
Gaussian joint distribution. 
For polynomials $f=f(t)\in \RR[t]$ homogeneous of degree $k$ we have
\begin{equation}\label{equation:FancyWickFormula}
\Ebold f(\tau)=\left\{\begin{array}{rl}
0&\mbox{if $k$ is odd,}\\
\displaystyle\frac{1}{m!}\left(\frac{1}{2}\sum_{i,j\in S}(\Ebold \tau_i\tau_j)\frac{\partial^2}{\partial t_i\partial t_j}\right)^mf(t)&\mbox{if $k=2m$ is even.}
\end{array}\right.
\end{equation}
\end{Lemma}
\proof If the family $\{\tau_i\}_{i\in S}$ is i.i.d. standard normal 
and $f(t)=(\sum_{i\in S}a_it_i)^k$ for a family of real constants $\{a_i\}_{i\in S}$ such that  $\sum_{i\in S}a_i^2=1$,
formula \eqref{equation:FancyWickFormula}
holds. Indeed, in that case the left side equals $\Ebold T^k$ for a standard normal random variable $T$,
and one can straightforwardly verify that the right side takes the same value.
But the polynomials $(\sum_{i\in S}a_it_i)^k$ span the subspace of $\RR[t]$ consisting of polynomials homogeneous of degree $k$, and
moreover, formula \eqref{equation:FancyWickFormula} is stable under homogeneous linear change of variable.  Thus formula \eqref{equation:FancyWickFormula} holds in general.
\qed

\subsubsection{The $\natural$-construction}
For any polynomial $g\in \RR[z]$ it is convenient 
to define a function $g^\natural:\Sym_n\rightarrow \RR$
which depends linearly on $g$ and which for $g$ homogeneous of degree $k$ is given by the formula
\begin{equation}\label{equation:NaturalDef}
g^\natural(Q)=
\left\{\begin{array}{rl}
0&\mbox{if $k$ is odd,}\\
\displaystyle\frac{1}{m!}\left(\frac{1}{2}\sum_{i,i'=1}^n\sum_{J=0}^{2N}Q(i,i')
(\Ebold \zeta_{i0}\zeta_{i'0})\frac{\partial^2}{\partial z_{iJ}\partial z_{i'J}}\right)^m
g(z)&\mbox{if $k=2m$ is even.}
\end{array}\right.
\end{equation}
It is clear that $g^\natural(Q)$ depends polynomially on the matrix entries of $Q$.

\begin{Lemma}\label{Lemma:HeatEquation}
For $g\in \RR[z]$ and $Q\in \Qfrak_n$ one has
\begin{equation}\label{equation:MomentFormula}
\Ebold g(\zeta\star Q)=g^\natural(Q).
\end{equation}
Furthermore, given also $\Gamma\subset\Bond_n$,
\begin{equation}\label{equation:HeatEquation}
\partial^\Gamma(g^\natural)=\left(\prod_{\{i,i'\}\in \Gamma}\Ebold \zeta_{i0}\zeta_{i'0}\right)(D^\Gamma g)^\natural.
\end{equation}
 \end{Lemma}
 \noindent Formula \eqref{equation:HeatEquation} is an algebraic variant of the {\em heat equation}.
 \proof Formula \eqref{equation:MomentFormula} follows immediately
 from Lemma \ref{Lemma:FancyWickFormula}. To prove \eqref{equation:HeatEquation} we may proceed
 by induction on $|\Gamma|$. The decisive case is then clearly that in which $\Gamma=\{e\}$
for some $e\in \Bond_n$.
In the latter special case differentiation on both sides of \eqref{equation:NaturalDef} immediately proves formula \eqref{equation:HeatEquation}.  \qed
\subsubsection{The independent copies trick}
We have
\begin{eqnarray}\label{equation:IndependentCopiesTrick}
\kappa\left(\left\{f_A(\zeta)\right\}_{A\in \Theta}\right)
&=&\sum_{\Pi\in [\Theta:\one_n]}\mu(\Pi:\one_n)
\prod_{B\in \Pi}
\Ebold \prod_{\begin{subarray}{c}
A\in \Theta\\
\mbox{\scriptsize s.t.}\;A\subset B
\end{subarray}}f_A(\zeta)
\\
\nonumber&=&\sum_{\Pi\in [\Theta:\one_n]}\mu(\Pi:\one_n)\prod_{B\in \Pi}
\Ebold \prod_{\begin{subarray}{c}
A\in \Theta\\
\mbox{\scriptsize s.t.}\;A\subset B
\end{subarray}}f_A(\zeta\star [\Pi])\\
\nonumber&=&\sum_{\Pi\in [\Theta:\one_n]}\mu(\Pi:\one_n)\Ebold f(\zeta\star [\Pi]).
\end{eqnarray}
The first step of the calculation is an application of formula
\eqref{equation:CumulantDefTer} and the remaining steps exploit
the covariance structure of the family $\zeta\star Q$ in a straightforward way.
The last step of the calculation is an instance of the commonly used
``independent copies trick'' whereby one writes a product of expectations
of random variables as the expectation of a product of independent copies of the variables.

\subsubsection{Application of the BKAR formula}
We have
\begin{eqnarray*}
(\mbox{LHS of \eqref{equation:MainTool}})
&=&\sum_{\Pi\in [\Theta:\one_n]}
\mu(\Pi:\one_n)\Ebold f(\zeta\star [\Pi])\;=\;\sum_{\Pi\in [\Theta:\one_n]}
\mu(\Pi:\one_n)f^\natural([\Pi])\\
\nonumber&=&\sum_{\Gamma\in \Tree_n(\Theta)}\int (\partial^\Gamma f^\natural)(Q)\PP^\Theta_\Gamma(dQ)\\
\nonumber&=&\sum_{\Gamma\in \Tree_n(\Theta)}\left(\prod_{\{i,i'\}\in \Gamma}(\Ebold \zeta_{i0}\zeta_{i'0})\right)\int (D^\Gamma f)^\natural(Q)\PP^\Theta_\Gamma(dQ)\\
&=&(\mbox{RHS of \eqref{equation:MainTool}}).
\end{eqnarray*}
The steps are justified as follows.
\begin{enumerate}
\item[Step 1.] Formula \eqref{equation:IndependentCopiesTrick} (independent copies trick).
\item[Step 2.] Formula \eqref{equation:MomentFormula} (Wick formula in terms of differential operators).
\item[Step 3.] Theorem \ref{Theorem:BKARarboreal} (BKAR formula).
\item[Step 4.] Formula \eqref{equation:HeatEquation} (heat equation).
\item[Step 5.] Formula \eqref{equation:MomentFormula} (Wick formula again).
\end{enumerate}
The proof of Theorem \ref{Theorem:MainTool} is complete. \qed

\subsection{Refinement of formula \eqref
{equation:PreMotzkinExpansion}}
\label{subsection:FieldObservation}
Under this heading we apply Theorem \ref{Theorem:MainTool}
to expand the right side of \eqref{equation:PreMotzkinExpansion} in a refined way.
\subsubsection{Opening the brackets}\label{subsubsection:OpenBrackets}
Fix a numerical partition $\lambda$ 
and a positive integer $N$.
Let $n=|\lambda|$ and $\ell=\ell(\lambda)$.
Fix $\theta\in S_n$ such that $\theta\sim \lambda$ 
and let $\Theta=\Orbit_n(\theta)\in \Part_n$. 
We then have
$$\Mfrak_{\lambda,N}
=\sum_{h:\langle n\rangle \rightarrow \langle N\rangle}
\kappa \left(\left\{\prod_{i\in A}\Tri_N(h(i),h(\theta(i)))\right\}_{A\in \Theta}\right)
$$
after opening the brackets in formula \eqref{equation:PreMotzkinExpansion} in evident fashion.
Let
\begin{equation}
\Motz_n^N(\theta)=\left\{
h:\langle n\rangle \rightarrow \langle N\rangle\bigg\vert
\max_{i\in \langle n\rangle}|h(\theta(i))-h(i)|\leq 1\right\}.
\end{equation}
\begin{figure}[hbt]
\centering
\includegraphics[scale=1]{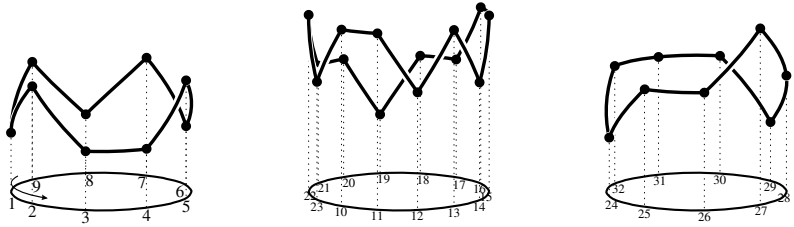}
\caption{
An $h\in {\rm{Motz}^{3}_{32}}(\theta)$ 
for a permutation $\theta\in S_{32}$ with three 
disjoint cycles.
\label{MotzN:Fig}
}
\end{figure}
Now a sequence of integers with increments in the set $\{0,\pm 1\}$ is often called a {\em Motzkin path}.
And it is evident that members of $\Motz_n^N(\theta)$ are  collections of closed Motzkin paths indexed by the blocks of the set partition $\Orbit_n(\theta)$. This is the rationale for our notation
$\Motz_n^N(\theta)$. Figure \ref{MotzN:Fig} provides an illustration.
 Given any function $h:\langle n\rangle\rightarrow\ZZ$, let
\begin{equation}
\label{equation:Jh}
J_h(\epsilon)=\{i\in \langle n\rangle\mid h(\theta(i))=h(i)+\epsilon\}\;\;
\mbox{for $\epsilon\in \{0,\pm 1\}$.}
\end{equation}
We then have
\begin{equation}\label{equation:MotzkinExpansion}
\Mfrak_{\lambda,N}
=\sum_{h\in \Motz_n^N(\theta)}\kappa\left(\left\{\prod_{i\in J_h(0)\cap A}\xi_{h(i)}\cdot \prod_{i\in  J_h(1)\cap A}\eta_{h(i)}\right\}_{A\in \Theta}\right)
\end{equation}
since 
$$\Tri_N(i+1,i)=0,\;\;\Tri_N(i,i)=\xi_i,\;\;\Tri_N(i,i+1)=\eta_i,$$
and otherwise for $|i-j|>1$ one has $\Tri_N(i,j)=0$.

\subsubsection{Specialization of Theorem \ref{Theorem:MainTool}}
Fix $h\in \Motz_n^N(\theta)$ arbitrarily. We will be considering not just one instance
of Theorem \ref{Theorem:MainTool} but rather a family of such instances indexed by $h$. For $A\in \Theta$ let
$$
f_A^h=\prod_{i\in J_h(0)\cap A}z_{i0}\cdot \prod_{i\in J_h(1)\cap A}\sum_{j=1}^{2h(i)}\frac{z_{ij}^2}{2}
\in\RR[\{\{z_{ij}\}_{i\in A}\}_{j=0}^{2N}].
$$
In turn let
\begin{equation}
f^h=\prod_{A\in \Theta}f_A^h\;=\;
\prod_{i\in J_h(0)}z_{i0}\cdot \prod_{i\in J_h(1)}\sum_{j=1}^{2h(i)}\frac{z_{ij}^2}{2}
\in \RR[\{\{z_{ij}\}_{i=1}^n\}_{j=0}^{2N}]=\RR[z].
\end{equation}
Let 
$$
\{\{\xi_{ij}\}_{i=1}^\infty\}_{j=0}^\infty
$$
be an i.i.d. family of standard normal random variables.
In turn consider the centered Gaussian family
\begin{equation}
\zeta^h=\left\{\left\{\zeta_{ij}^h\right\}_{i=1}^n\right\}_{j=0}^{2 N}=\left\{\left\{\xi_{h(i),j}\right\}_{i=1}^n\right\}_{j=0}^{2 N}.
\end{equation}
Note that by definition of $\zeta^h$ we have
\begin{equation}\label{equation:These1}
\Ebold \zeta_{i_1j_1}^h\zeta_{i_2j_2}^h=\Ebold \xi_{h(i_1),j_1}\xi_{h(i_2),j_2}=\delta_{h(i_1),h(i_2)}\delta_{j_1j_2}.
\end{equation}
For each $Q\in \Qfrak_n$, let $\zeta^h\star Q$ denote the family arising from $\zeta^h$ by $Q$-recoupling. 

Here is the promised refined expansion of the right side of \eqref{equation:PreMotzkinExpansion}.
\begin{Proposition}\label{Proposition:ToolCrush} 
Notation and assumptions are as above. We have
\begin{eqnarray}\label{equation:MotzkinExpansionSteroid}
\Mfrak_{\lambda,N}
&=&\sum_{(\Gamma,h)\in \Tree_n(\Theta)\times \Motz_n^N(\theta)}\\
\nonumber&&\left(\prod_{\{i_1,i_2\}\in \Gamma}\delta_{h(i_1),h(i_2)}\right)
\int \Ebold\left[(D^\Gamma f^h)(\zeta^h\star Q)\right]\, \PP^\Theta_\Gamma(\dd Q).
\end{eqnarray}
Furthermore, the summand on the right side of \eqref{equation:MotzkinExpansionSteroid} 
indexed by $(\Gamma,h)$ vanishes unless the following four conditions hold:
\begin{eqnarray}
\label{equation:MotzScrambleRedux1}
&&\mbox{$h$ is constant on each  $e\in \Gamma$.}\\
\label{equation:MotzScrambleRedux2}
&&
\mbox{$e\subset J_h(0)$ or $e\subset J_h(1)$ for each $e\in \Gamma$.}\\
\label{equation:MotzScrambleRedux3}
&&\mbox{For each $i\in J_h(0)$ there exists at most one $e\in \Gamma$ such that $i\in e$.}\\
\label{equation:MotzScrambleRedux4}
&&\mbox{For each $i\in \langle n\rangle$ there exist at most two $e\in \Gamma$ such that $i\in e$.}
\end{eqnarray}\end{Proposition}
\proof We begin by proving formula \eqref{equation:MotzkinExpansionSteroid}.
To do so it suffices by \eqref{equation:MotzkinExpansion} to fix \linebreak $h\in \Motz_n^N(\theta)$ arbitrarily
and to prove that
\begin{eqnarray}\label{equation:MotzkinExpansionSteroidPrime}
&&\kappa\left(\left\{\prod_{i\in J_h(0)\cap A}\xi_{h(i)}\cdot \prod_{i\in J_h(1)\cap A}\eta_{h(i)}\right\}_{A\in \Theta}\right)\\
\nonumber &=&\sum_{\Gamma\in \Tree_n(\Theta)}\left(\prod_{\{i_1,i_2\}\in \Gamma}\delta_{h(i_1),h(i_2)}\right)\int \Ebold\left[(D^\Gamma f^h)(\zeta^h\star Q)\right]\, \PP^\Theta_\Gamma(\dd Q).
\end{eqnarray}
Now we are free to replace the family  $\{\xi_i\}_{i=1}^N\cup\{\eta_i\}_{i=1}^{N-1}$ 
appearing in the definition of the tridiagonal matrix $\Tri_N$ by any other family with the same joint law.
Thus we may assume without loss of generality that 
\begin{equation}\label{equation:These2}
\xi_i=\xi_{i0}\;\;\mbox{and}\;\;\eta_i=\sum_{j=1}^{2i} \frac{\xi_{ij}^2}{2}.
\end{equation}
Relations \eqref{equation:These1} and \eqref{equation:These2} 
taken into account, it is clear that \eqref{equation:MotzkinExpansionSteroidPrime}
is a specialization of Theorem \ref{Theorem:MainTool} and thus holds.
Thus in turn formula \eqref{equation:MotzkinExpansionSteroid} indeed holds.

Now fix a pair $(\Gamma,h)$ such that \eqref{equation:MotzScrambleRedux1} holds
and moreover $D^\Gamma f^h\neq 0$. It will be enough to show that for this pair $(\Gamma,h)$ statements \eqref{equation:MotzScrambleRedux2}---\eqref{equation:MotzScrambleRedux4} hold.
By opening the brackets 
in the definition of $f^h$ we infer the existence of a function $g:J_h(1)\rightarrow \ZZ$
such that $0<g(i)\leq 2h(i)$ for $i\in J_h(1)$ and such that the monomial
$$Z=\prod_{i\in J_h(0)}z_{i0}\cdot \prod_{i\in J_h(1)}z_{i,g(i)}^2$$
satisfies $D^\Gamma Z\neq 0$. 
Failure of \eqref{equation:MotzScrambleRedux2} would entail existence of $e\in \Gamma$
such that either $e$ meets both $J_h(0)$ and $J_h(1)$
or else $e$ meets $J_h(-1)$.
In both cases, in the former by \eqref{equation:FirstCull}
and in the latter by \eqref{equation:SecondCull},
we would have $D^\Gamma Z=0$, which is a contradiction.
Thus \eqref{equation:MotzScrambleRedux2} holds.
Failure of \eqref{equation:MotzScrambleRedux3} would entail existence of distinct $e,e'\in \Gamma$
such that $e\cap e'\cap J_h(0)\neq \emptyset$. In this case we would have $D^\Gamma Z=0$
by \eqref{equation:SecondCull}, which is again a contradiction. Thus \eqref{equation:MotzScrambleRedux3} holds.
Failure of \eqref{equation:MotzScrambleRedux4} would entail existence of distinct $e_1,e_2,e_3\in \Gamma$
such that $e_1\cap e_2\cap e_3 \neq \emptyset$, in which case $D^\Gamma Z=0$
by \eqref{equation:SecondCull},  which is yet again a contradiction. Thus \eqref{equation:MotzScrambleRedux4} holds.
The proof of Proposition \ref{Proposition:ToolCrush} is complete.
\qed

\section{Linear forests, cycle-cut permutations and Goulden-Jackson pairs}
\label{section:GraphPreparation}
In this section for the sake of clarity we hold ourselves somewhat aloof from the proof of Theorem \ref{Theorem:MainResult} and develop some simple concepts on their own terms.  All these concepts
are motivated by Proposition \ref{Proposition:ToolCrush} and they will be deployed in \S\ref{section:TridiagonalCumulantCalc} below to clinch the proof of Theorem \ref{Theorem:MainResult}.

\subsection{Linear forests}
The notion developed under this heading is directly
motivated by statement \eqref{equation:MotzScrambleRedux4} of 
Proposition \ref{Proposition:ToolCrush} above.

\subsubsection{Definition} Let $\Gamma\subset \Bond_n$ be a subset.
We call $\Gamma$ a {\em linear forest}
if the graph $\Gfrak(\zero_n,\Gamma)$ is circuitless (i.e., a forest) and every vertex of the graph $\Gfrak(\zero_n,\Gamma)$ has degree at most $2$.

\subsubsection{The set $\Tree_n^\LF(\Theta)$}
For $\Theta\in \Part_n$, 
let $\Tree^{\lf}_n(\Theta)$ 
denote the subset of $\Tree_n(\Theta)$ consisting of linear forests. 
Now suppose $\Gamma\in \Tree_n(\Theta)$ is given.
Then $\Gamma$ necessarily has the property that the graph $\Gfrak(\zero_n,\Gamma)$ is circuitless.
Thus for $\Gamma\in \Tree_n(\Theta)$,
one has $\Gamma\in \Tree_n^\LF(\Theta)$ if and only if every 
vertex of $\Gfrak(\zero_n,\Gamma)$ has degree at most $2$.

\subsubsection{Decomposition of linear forests into connected components}
We call a linear forest $\Gamma\subset \Bond_n$ {\em connected}
if the graph $\Gfrak(\zero_n,\Gamma)$ has exactly one connected component not reducing to an isolated vertex.
Every connected linear forest $\Gamma\subset \Bond_n$ is of the form
\begin{equation}\label{equation:ConnectedExamples}
\Gamma=\{\{i_1,i_2\},\dots,\{i_{m-1},i_m\}\}\;\mbox{for $m\geq 2$
and distinct $i_1,\dots,i_m\in \langle n\rangle$.}
\end{equation}
Note furthermore that the sequence $i_1,\dots,i_m$ is uniquely determined by $\Gamma$ up to a reversal
of the order of the sequence.
It is clear that every linear forest $\Gamma$ has a disjoint union decomposition
$$\Gamma=\Gamma_1\cup\cdots \cup\Gamma_k$$
unique up to ordering of the sets in the decomposition,
where each set $\Gamma_i$ is a connected linear forest
and the union  
$$\cup \Gamma=(\cup\Gamma_1)\cup\cdots\cup(\cup\Gamma_k)$$ is also disjoint.
Each set $\Gamma_i$ is called a {\em connected component} of $\Gamma$.
We say that $\Gamma=\bigcup\Gamma_i$ is the {\em decomposition} of $\Gamma$
into its connected components.

\subsubsection{The  boundary of a linear forest}
Let $\Gamma\subset \Bond_n$ be a linear forest.
We define the {\em boundary} $\partial\Gamma$ to be the set of
unordered pairs of the form $\{a,b\}$ where $\{a\}$ and $\{b\}$ are distinct degree one vertices 
of the forest $\Gfrak(\zero_n,\Gamma)$ joined by some walk.  The members 
of $\partial \Gamma$ are in evident bijection with the connected components of $\Gamma$.

\begin{figure}[hbt]
\centering
\includegraphics[scale=1]{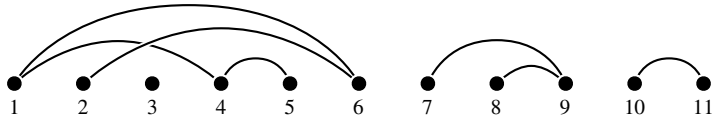}
\caption{
The graph $\mathfrak G(\mathbf{0}_{11},\Gamma)$ 
is depicted above for
the linear forest ${\Gamma=\{\{1,4\},\{1,6\},\{2,6\},\{4,5\},\{7,9\},\{8,9\},\{10,11\}\}}\subset \Bond_{11}$. 
The connected components of $\Gamma$ are  
$\Gamma_1=\{\{1,4\},\{1,6\},\{2,6\},\{4,5\}\}$, 
$\Gamma_2=\{\{7,9\},\{8,9\}\}$, and 
$\Gamma_3=\{\{10,11\}\}$.  
The boundary of $\Gamma$ is $\partial \Gamma=\{\{2,5\},\{7,8\},\{10,11\}\}$.
\label{LF:Fig}
}
\end{figure}

\subsection{Cycle-cut permutations}
We next introduce a notion which is nearly equivalent to that of a linear forest.
\subsubsection{Cycle-cuttings}\label{subsubsection:CycleCut}
 Let $\sigma\in S_n$ be any permutation.
A subset of $\supp \,\sigma$ 
intersecting each  $\sigma$-orbit contained in $\supp\,\sigma$ in exactly one point
will be called a {\em cycle-cutting}. A pair $(\sigma,A)$ consisting of $\sigma\in S_n$
and a cycle-cutting $A\subset \langle n\rangle$ of $\sigma$ will be called a {\em cycle-cut permutation}.
 Now let $\lambda\vdash n$ index the conjugacy class of $\sigma$.
We define
\begin{equation}
\mbold(\sigma)=\mbold(\lambda)= \prod_i i^{m_i(\lambda)}=\prod_{i=1}^{\ell(\lambda)}\lambda_i.
\end{equation}
Note that $\sigma$ has exactly $\mbold(\sigma)$ cycle-cuttings.
Given also $i\in \langle n\rangle$, let $\mbold(\sigma,i)$ denote the cardinality of the $\sigma$-orbit
to which $i$ belongs.  Note that 
\begin{equation}\label{equation:Embolden}
\mbold(\sigma)=\prod_{a\in A}\mbold(\sigma,a)
\end{equation}
for any cycle-cutting $A$ of $\sigma$.
\subsubsection{Construction of linear forests from cycle-cut permutations}
Given a cycle-cut permutation $(\sigma,A)$ of $\langle n\rangle$, let 
$$\LF(\sigma,A)=\{\{i,\sigma(i)\}\mid i\in (\supp\, \sigma)\setminus A\}\subset \Bond_n,$$
which is clearly a linear forest. 
For each cycle-cut permutation $(\sigma,A)$ of $\langle n\rangle$
and associated linear forest $\Gamma=\LF(\sigma,A)$
it is furthermore clear  that 
\begin{eqnarray}
&&\cup \Gamma=\supp \,\sigma,\;
|\Gamma|=n-\ell(\sigma),\;
\partial \Gamma=\{\{a,\sigma(a)\}\mid a\in A\},\;\mbox{and}\;\\
&&|\partial \Gamma|=|A|=\ell(\sigma)-|\{i\in \langle n\rangle\mid \sigma(i)=i\}|\\
\nonumber&=&\mbox{the number of connected components of $\Gamma$.}
\end{eqnarray}

\begin{figure}[hbt]
\centering
\includegraphics[scale=1]{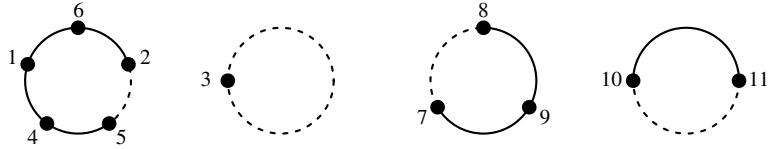}
\caption{
For the cycle-cut permutation 
$(\sigma,A)$ depicted above
 where $\sigma=(1,4,5,2,6)(7,9,8)(10,11)\in S_{11}$ and 
$A=\{5,8,10\}$,
the linear forest $\Gamma=\LF(\sigma,A)$  is ${\Gamma=\{\{1,4\},\{1,6\},\{2,6\},\{4,5\},\{7,9\},\{8,9\},\{10,11\}\}}$ and  
the boundary is $\partial \Gamma=\{\{2,5\},\{7,8\},\{10,11\}\}$.
\label{CycleCutting:Fig}
}
\end{figure}

\begin{figure}[hbt]
\centering
\includegraphics[scale=1]{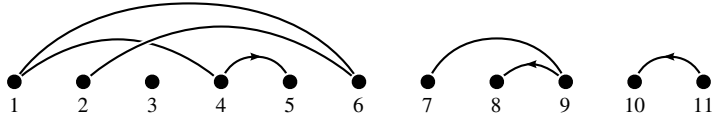}
\caption{The oriented linear forest associated with the permutation
$\sigma=(1,4,5,2,6)(7,9,8)(10,11)\in S_{11}$ and the cycle-cutting
$A=\{5,8,10\}$.
\label{LFChoice:Fig}
}
\end{figure}

The notions of cycle-cut permutation and of linear forest
are equivalent up to some manageable powers of $2$, as the next  lemma explains.
\begin{Lemma} \label{Lemma:BreakCycle}
For each linear forest $\Gamma\subset \Bond_n$, the set of cycle-cut permutations $(\sigma,A)$
of $\langle n\rangle$ such that $\Gamma=\LF(\sigma,A)$ has cardinality  $2^{k}$
where $k$ is the number of connected components of $\Gamma$ (and hence $k=|A|$).
 \end{Lemma}
\proof Given a linear forest $\Gamma$,
let us (temporarily, just within this proof) call a choice of point from each member of the boundary $\partial \Gamma$
an {\em orientation}. And in turn (again, temporarily) let us call $\Gamma$ an {\em oriented}
linear forest if it is equipped with an orientation.
The notion of cycle-cut permutation
is precisely equivalent to the notion of oriented linear forest,
with cycle-cuttings corresponding one-to-one with orientations.
See Figure \ref{LFChoice:Fig} for an illustration.
Obviously $\Gamma$ has exactly $2^{|\partial\Gamma|}$ orientations.
Thus the lemma holds.  (Going forward we will not make further use of oriented linear forests.) \qed

\subsection{Relationship between linear forests 
and Goulden-Jackson pairs}
We now make the decisive linkage between on the one hand combinatorial objects related to the BKAR formula and on the other hand Goulden-Jackson pairs.

\begin{Proposition}\label{Proposition:LoopCuttingGJ}  
Fix $\theta\in S_n$ and let $\Theta=\Orbit_n(\theta)$.
Fix a cycle-cut permutation $(\sigma,A)$ of $\langle n\rangle$
and let $\Gamma=\LF(\sigma, A)\subset \Bond_n$.
Then $\sigma\in\GJ_n(\theta)$ iff
$\Gamma\in \Tree^\lf_n(\Theta)$.
\end{Proposition}
\proof Consider the {\em Cayley graph}
$$\Cfrak=\Gfrak(\zero_n,\{\{i,\theta(i)\},\{i,\sigma(i)\}\mid i\in \langle n\rangle\}).$$
Then the subgroup of $S_n$ generated by $\theta$ and $\sigma$ acts transitively on $\langle n\rangle$
iff $\Cfrak$ is connected. It is easy to check in turn that $\Cfrak$ is connected iff $\Gfrak(\Theta,\Gamma)$ is connected.

Suppose now that we have $\sigma\in \GJ_n(\theta)$, i.e., $(\theta,\sigma)\in \GJ_n$.
Then  
$$|\Gamma|=n-\ell(\sigma)=\ell(\theta)-1=|\Theta|-1.$$
Furthermore $\Cfrak$ is connected since $\ell(\sigma\theta)=1$
and hence $\Gfrak(\Theta,\Gamma)$ is connected. 
Thus $\Gfrak(\Theta,\Gamma)$ is a tree and hence $\Gamma\in \Tree^\lf(\Theta)$.

Suppose now rather that $\Gamma\in \Tree^\lf(\Theta)$ and hence that $\Gfrak(\Theta,\Gamma)$
is a tree. Then  we have 
$$n-\ell(\sigma)=|\Gamma|=|\Theta|-1=\ell(\theta)-1$$
and hence $\ell(\sigma)+\ell(\theta)=n+1$ holds. It remains only verify $\ell(\sigma\theta)=1$.
In any case, $\Gfrak(\Theta,\Gamma)$ is connected, hence
$\Cfrak$ is connected and hence $\sigma$ and $\theta$ generate a subgroup
of $S_n$ acting transitively on $\langle n\rangle$. 
Lemma \ref{Lemma:RiemannHurwitz}  immediately below then yields the bound $\ell(\sigma\theta)\leq 1$,
which finishes the proof.
\qed

\begin{Lemma}\label{Lemma:RiemannHurwitz}
Let $\sigma,\theta\in S_n$ be permutations together generating a subgroup of $S_n$ acting
transitively on $\langle n\rangle$. Then $\ell(\sigma)+\ell(\theta)+\ell(\sigma\theta)\leq n+2$.
\end{Lemma}
\proof The lemma reiterates \cite[Thm. 3.6, p. 421]{CoriMachi} in different notation.
For readers familiar with the theory of compact Riemann surfaces, we supplement
this reference with the following brief explanation.
From the permutations $\theta$ and $\sigma$ one knows how to construct a compact Riemann surface
of genus $g$ ($=$ number of handles) presented as an $n$-sheeted covering of the Riemann sphere
branched only at $0$, $1$ and $\infty$ such that the {\em Riemann-Hurwitz formula}
$$2g-2=-2n+(n-\ell(\theta))+(n-\ell(\sigma))+(n-\ell(\theta\sigma))$$
holds. The desired inequality follows simply from the fact that $g\geq 0$.
\qed

The preceding theory provides valuable information
about integrals against the measure $\PP^\Theta_\Gamma$ 
of certain simple functions.
\begin{Proposition}\label{Proposition:LastCut}
Fix $(\theta,\sigma)\in\GJ_n$.
Let $\Theta=\Orbit_n(\theta)$.
Let $A$ be a cycle-cutting for $\sigma$.  
Let $\Gamma=\LF(\sigma,A)\in \Tree^{\LF}_n(\Theta)$. 
For $i,j\in \langle n\rangle$ let $\Gamma(i,j)\subset \Gamma$ be as defined in Lemma \ref{Lemma:PTreeChar}.
Let $X\in \Qfrak_n$ be a random matrix with law $\PP_\Gamma^\Theta$.
(i) The set $\Gamma$ is the disjoint union of the sets $\Gamma(a,\sigma(a))$ for $a\in A$.
(ii) The family $\{X(a,\sigma(a))\}_{a\in A}$ of random variables is independent.
(iii) $\Ebold X(a,\sigma(a))=\frac{1}{\mbold(\sigma,a)}$ for $a\in A$. 
(iv) One has 
\begin{equation}\label{equation:QInt}
\Ebold \prod_{b\in B}X(b,\sigma(b))=1\bigg/\prod_{b\in B}\mbold(\sigma,b)
\end{equation}
for any subset $B\subset A$.
\end{Proposition}
\proof 
 Let $\Gamma=\bigcup_{a\in A} \Gamma_a$ be the unique decomposition of $\Gamma$ into its connected components labeled so that $\partial\Gamma_a=\{a,\sigma(a)\}$ for $a\in A$.
It is not hard to see that for each $a\in A$ one has
 $\Gamma_a=\Gamma(a,\sigma(a))$. Thus statement (i) holds.
Statement (ii)  follows via Lemma \ref{Lemma:PTreeChar} from statement (i).
Statement (iii) follows from Lemma \ref{Lemma:PTreeChar}
and the undergraduate-level remark that 
for random variables $U_1,\dots,U_k$ i.i.d. uniform in $(0,1)$
one has $\Ebold \min_{i=1}^kU_i=\frac{1}{k+1}$. Statement (iv) follows immediately
from statements (ii) and (iii).
\qed

\begin{figure}[hbt]
\centering
\includegraphics[scale=1]{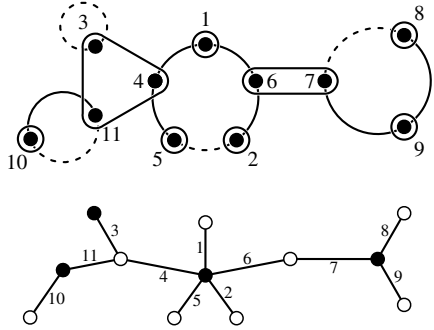}
\caption{
The top figure shows the cycle-cut permutation 
$(\sigma,A)$ and the set partition 
$\Theta={\rm Orb}_n(\theta)$ for 
$\theta=(3,11,4)(6,7)$, $\sigma=(1,4,5,2,6)(7,9,8)(10,11)$ in $S_{11}$, and 
$A=\{5,8,10\}$.  
The corresponding Shabat-Voevodsky tree for $(\theta,\sigma)$ is the bottom figure.
}
\end{figure}

\subsection{Objects related to $\dMotz_n(\theta,\sigma)$}\label{subsection:Alternative}
Having developed above an interpretation of statement \eqref{equation:MotzScrambleRedux4}
of Proposition \ref{Proposition:ToolCrush} in group-theoretical terms, we turn next to the task
of providing an analogous interpretation of
 statements \eqref{equation:MotzScrambleRedux1}---\eqref{equation:MotzScrambleRedux3}.

\subsubsection{The set $\Motz_n(\theta,\sigma)$}
For $(\theta,\sigma)\in \GJ_n$, 
let  
\begin{equation}
\Motz_n(\theta,\sigma)=
\{h:\langle n\rangle\rightarrow\ZZ\mid
h\circ \sigma=h\;\mbox{and}\;
h\circ \theta-h\in \dMotz_n(\theta,\sigma)\}.
\end{equation}
In this setting we think of the function $h\circ \theta-h$ as a sort of derivative of $h$.
From this definition one immediately deduces
the following statement:
\begin{equation}
\label{equation:MotzStability}
h\in \Motz_n(\theta,\sigma)\Leftrightarrow h+c\in \Motz_n(\theta,\sigma)\;\mbox{for constants $c\in \ZZ$.}
\end{equation}
For each positive integer $N$ we also define
\begin{equation}\label{equation:MotzupperNotation}
\Motz_n^N(\theta,\sigma)=\Motz_n(\theta,\sigma)\cap \Motz_n^N(\theta).
\end{equation}

\subsubsection{``Tilde versions'' of the preceding definitions}
Let $(\theta,\sigma)\in \GJ_n$. 
Let
$$\widetilde{\dMotz}_n(\theta,\sigma)\supset \dMotz_n(\theta,\sigma)$$
be the superset consisting of $g:\langle n\rangle\rightarrow\ZZ$ satisfying \eqref{equation:Interaction0}--\eqref{equation:Interaction3.5}
but perhaps not satisfying \eqref{equation:Interaction4}.
In turn, let 
\begin{equation}
\widetilde{\Motz}_n(\theta,\sigma)=
\{h:\langle n\rangle\rightarrow\ZZ\mid
h\circ \sigma=h\;\mbox{and}\;
h\circ \theta-h\in \widetilde{\dMotz}_n(\theta,\sigma)\}.
\end{equation}
Note that the variant of \eqref{equation:MotzStability} with $\widetilde{\Motz}_n(\theta,\sigma)$
in place of $\Motz_n(\theta,\sigma)$ still holds.
Note the trivial but important relation
\begin{equation}\label{equation:UltraTrivial}
\Motz_n(\theta,\sigma)=\{h\in \widetilde{\Motz}_n(\theta,\sigma)\mid J_h(0)\subset \supp\,\sigma\}
\end{equation}
where $J_h(\epsilon)$ is as defined on line \eqref{equation:Jh}.
We also define
\begin{equation}
\widetilde{\Motz}_n^N(\theta,\sigma)=\Motz_n^N(\theta)\cap\widetilde{\Motz}_n(\theta,\sigma).
\end{equation}

\begin{Proposition}\label{Proposition:MotzScrambleInterp}
Fix $(\theta,\sigma)\in \GJ_n$, a cycle-cutting $A$ of $\sigma$
and $h\in \Motz_n^N(\theta)$. Let $\Gamma=\LF(\sigma,A)$.
If the pair $(\Gamma,h)$ satisfies statements \eqref{equation:MotzScrambleRedux1}---\eqref{equation:MotzScrambleRedux3}, then \linebreak $h\in \widetilde{\Motz}_n^N(\theta,\sigma)$.
\end{Proposition}
\proof 
Statement \eqref{equation:MotzScrambleRedux1} implies $h\circ \sigma=h$.
Let $g=h\circ \theta-h$. It remains only to show that
$g\in \widetilde{\dMotz}_n^N(\theta)$.
The definition of $\Motz_n^N(\theta)$ implies that
$g$ satisfies \eqref{equation:Interaction0}.
Clearly, $g$ satisfies \eqref{equation:Interaction1}.
 Statement  \eqref{equation:MotzScrambleRedux2} implies
 that $g$ satisfies \eqref{equation:Interaction2} and \eqref{equation:Interaction3}.
Statement \eqref{equation:MotzScrambleRedux3} implies that $g$ satisfies \eqref{equation:Interaction3.5}.
 Thus we indeed have $g\in \widetilde{\dMotz}_n^N(\theta)$ and hence
$h\in \widetilde{\Motz}_n^N(\theta,\sigma)$.  
\qed

\subsection{``Integration'' on Goulden-Jackson pairs}
We explain the sense in which each element of $\dMotz(\theta,\sigma)$
has an antiderivative. (See Proposition \ref{Proposition:EquiNumerous} below.)

\begin{Lemma}\label{Lemma:StupidBound}   Fix $(\theta,\sigma)\in \GJ_n$
and a function
$h:\langle n\rangle\rightarrow \ZZ$ 
such that $h\circ \sigma=h$.
Then we have
\begin{equation}
\max_{i,j\in \langle n\rangle}|h(i)-h(j)|\leq (n-1)\max_{i\in \langle n\rangle}|h(\theta(i))-h(i)|.
\end{equation}
In particular, if $h\circ \theta-h=0$, then $h$ is constant.
\end{Lemma}
\proof By hypothesis 
$|h\circ(\sigma\theta)-h|=|h\circ \theta-h|$ and $\ell(\sigma\theta)=1$,
whence the bound.
\qed

\begin{Lemma}\label{Lemma:hIntegration}
Fix $(\theta,\sigma)\in \GJ_n$. Let
$g:\langle n\rangle\rightarrow \ZZ$ be a function averaging 
to $0$ on each $\theta$-orbit.
Then there exists a function
$h:\langle n\rangle\rightarrow \ZZ$ 
such that $h\circ \sigma=h$ and $h\circ \theta-h=g$.
\end{Lemma}
\proof Consider again the Cayley graph 
$$\Cfrak=\Gfrak(\zero_n,\{\{i,\theta(i)\},\{i,\sigma(i)\}\mid i\in \langle n\rangle\})$$
introduced in the proof of Proposition \ref{Proposition:LoopCuttingGJ}.
Since $\ell(\theta\sigma)=1$, it is clear that $\Cfrak$ is connected.
Let $A$ (resp., $B$) be a cycle-cutting for $\theta$ (resp., $\sigma$).
Let
$$\Tfrak=\Gfrak(\zero_n,\{\{i,\theta(i)\}\mid i\in (\supp\,\theta)\setminus A\}
\cup\{\{j,\sigma(j)\}\mid j\in (\supp\,\sigma)\setminus B\}).$$
It is clear that any two distinct vertices of $\Tfrak$ joined by a walk in $\Cfrak$ remain joined by some walk in $\Tfrak$. (For every bridge knocked out an alternate route has been preserved.)
Thus $\Tfrak$ is connected.
Furthermore, $\Tfrak$ has no more than $n-1$ edges because
$$|(\supp\,\theta)\setminus A|+|(\supp\,\sigma)\setminus B|=(n-\ell(\theta))+(n-\ell(\sigma))=n-1.$$
Thus $\Tfrak$ is a tree spanning $\Cfrak$. In particular, $\Tfrak$ has exactly $n-1$ edges.
Now (so to speak) every vector field on a tree is the gradient of a potential and this statement
holds over $\ZZ$.
Thus there exists some function $h:\langle n\rangle\rightarrow \ZZ$ such that 
$$h(\theta(i))-h(i)=g(i)\;\;
\mbox{for $i\in A$ and}\;\;h(\sigma(j))-h(j)=0\;\;\mbox{for $j\in B$.}$$
Clearly, we have $h\circ \sigma=h$.
Finally,  since $g$ averages to zero on each $\theta$-orbit, we have $h\circ \theta-h=g$. Thus $h$ exists.
\qed

\begin{figure}[hbt]
\centering
\includegraphics[scale=1]{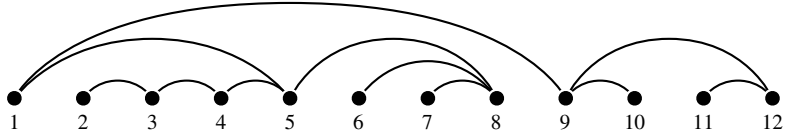}
\caption{
The tree $\mathfrak T$ from Lemma 3.5.4 
for the cycle-cut permutations $(\theta,A)$ and $(\sigma,B)$
where 
$(\theta,\sigma)\in \GJ_{12}$,
$\theta=(1,9)(2,3,4,5)(6,7,8)$, 
$\sigma=(1,5,8)(9,10,11,12)$,
$A=\{5,6,9\}$, and 
$B=\{8,10\}$.  
\label{TwoCycleCutTree:Fig}
}
\end{figure}

\begin{Proposition}\label{Proposition:EquiNumerous}
For $(\theta,\sigma)\in \GJ_n$ and  $g\in \dMotz_n(\theta,\sigma)$
there exists unique $h\in \Motz_n(\theta,\sigma)$
such that $h\circ \theta-h=g$, $h\circ \sigma=h$ and $h(1)=0$.
\end{Proposition}
\noindent Thinking of $h\circ \theta-h$ as the derivative of $h$,
this finally is the rationale for the peculiar  notation $\dMotz(\theta,\sigma)$.
\proof Lemma \ref{Lemma:StupidBound},  Lemma \ref{Lemma:hIntegration}, 
and statement \eqref{equation:MotzStability} prove this. \qed

Finally, we record an easy consequence of Lemma \ref{Lemma:StupidBound} 
for  convenient reference.
\begin{Lemma}
For $(\theta,\sigma)\in \GJ_n$ we have crude bounds
\begin{equation}\label{equation:CrudeMotzBound}
\left|\left\{h\in \widetilde{\Motz}_n(\theta,\sigma)\bigg\vert h(1)=0\right\}\right|\leq 3^n\;\;
\mbox{and}\;\;
\left|\widetilde{\Motz}_n^N(\theta,\sigma)\right|\leq 3^nN.
\end{equation}
\end{Lemma}
\proof Lemma \ref{Lemma:StupidBound} implies the first bound.
The latter in turn implies the second bound. \qed

\subsection{A limit calculation}
Our main result under this heading explains the denominator
on the right side of  formula \eqref{equation:MainResult}.
(See  Proposition \ref{Proposition:IntMotz} below.)
\begin{Lemma}\label{Lemma:BalanceBooks}
For $(\theta,\sigma)\in \GJ_n$, a cycle-cutting $A$ of $\sigma$
and $h\in \widetilde{\Motz}_n(\theta,\sigma)$ we have
\begin{eqnarray}
\label{equation:SpinUp}
\frac{1}{\mbold(\sigma)}&=&2^{-|A|}\prod_{a\in J_h(1)\cap A}\frac{2}{\mbold(\sigma,a)}\;\;\mbox{and}\\
\label{equation:BalanceBooks}
\frac{n}{2}-\ell(\theta)+1&=&\frac{|J_h(0)\setminus \supp\,\sigma|}{2}+|J_h(1)\setminus \supp\, \sigma|
+|J_h(1)\cap A| \geq 0
\end{eqnarray}
\end{Lemma}
\proof 
 By \eqref{equation:Interaction0} and \eqref{equation:Interaction1} we evidently have
$$\frac{n}{2}=\frac{|J_h(0)|+|J_h(1)|+|J_h(-1)|}{2}=\frac{|J_h(0)|}{2}+|J_h(1)|.$$
By \eqref{equation:Interaction2} each of the sets $J_h(-1)$, $J_h(0)$ and $J_h(1)$  is $\sigma$-stable,
i.e., each is a union of $\sigma$-orbits.
 To abbreviate notation let $\Sigma=\supp \,\sigma$.
We have
$$
\ell(\theta)-1=n-\ell(\sigma)=|\Sigma \setminus A|
=|J_h(0)\cap (\Sigma \setminus A)|+|J_h(1)\cap (\Sigma \setminus A)|,
$$
at the first step by definition of a Goulden-Jackson pair,
at the second step as a consequence of the definition of a cycle-cutting
and at the last step by \eqref{equation:Interaction3}.
By \eqref{equation:Interaction3.5} we have $\mbold(\sigma,a)=2$ for $a\in A\cap J_h(0)$,
whence \eqref{equation:SpinUp} via \eqref{equation:Embolden}.
And furthermore, we have
$$
|J_h(0)\cap A|=\frac{|J_h(0)\cap\Sigma|}{2}=|J_h(0)\cap(\Sigma\setminus A)|.
$$ 
Formula \eqref{equation:BalanceBooks}
 can then be obtained by combining the three
displayed lines above.
\qed

\begin{Lemma}\label{Lemma:BalanceBooksBis}
For $\theta\in S_n$ such that
$\frac{n}{2}-\ell(\theta)+2\leq 0$ the set $\GJdM_n(\theta)$ is empty.
\end{Lemma}
\noindent In other words, statement \eqref{equation:AintZeroNeither} above holds.
\proof Supposing $\GJdM_n(\theta)$ is not empty,
there exists some $\sigma\in \GJ_n(\theta)$ and some $g\in \dMotz_n(\theta,\sigma)$.
In turn, by Proposition \ref{Proposition:EquiNumerous} there exists some $h\in \Motz_n(\theta,\sigma)$ such that $h\circ \theta-h=g$.
By Lemma \ref{Lemma:BalanceBooks} we would then have
$\frac{n}{2}-\ell(\theta)+2> 0$, which is a contradiction.
\qed

\begin{Proposition}\label{Proposition:IntMotz}
Let $(\theta,\sigma)\in \GJ_n$. Let $\ell=\ell(\theta)$.
Let $A$ be a cycle-cutting of $\sigma$. 
Let $N$ be a positive integer. Let $g\in \dMotz_n(\theta,\sigma)$.
Then we have
\begin{equation}\label{equation:FinalIntegration}
\left|
\frac{N^{\frac{n}{2}-\ell+2}}{\frac{n}{2}-\ell+2}-\sum_{
\begin{subarray}{c}
h\in \Motz_n^N(\theta,\sigma)\\
\textup{\mbox{\scriptsize s.t.}}\,h\circ \theta-h=g\;\;\;\;
\end{subarray}}\;\prod_{i\in J_h(1)
\setminus((\supp\,\sigma)\setminus A)}h(i)
\right|\\
\leq cN^{\frac{n}{2}-\ell+1},
\end{equation}
where the constant $c$ depends only on $n$.
\end{Proposition}
\proof
Let
$$
H=\sum_{
\begin{subarray}{c}
h\in \Motz_n^N(\theta,\sigma)\\
\textup{\mbox{\scriptsize s.t.}}\,h\circ \theta -h=g\;\;\;\;
\end{subarray}}\frac{1}{N}\prod_{i\in J_h(1)\setminus ((\supp\,\sigma)\setminus A)}\left(\frac{h(i)}{N}\right)\;\;
\mbox{and}\;\;
\nu=\frac{n}{2}-\ell+1.
$$
It will be enough to prove that
\begin{equation}\label{equation:FinalIntegrationBis}
\left|H-\int_0^1 t^{\nu}\,\dd t\right|\leq \frac{(n+1)^2}{N}.
\end{equation}
By  \eqref{equation:UltraTrivial} and \eqref{equation:BalanceBooks} we have
$$h\in \Motz_n^N(\theta,\sigma)\Rightarrow \nu=|J_h(1)\setminus \supp\,\sigma|+
|J_h(1)\cap A|=|J_h(1)\setminus((\supp\,\sigma) \setminus A)|.$$
 By  Lemma \ref{Lemma:StupidBound} and the definitions we have 
$$h\in \Motz_n(\theta,\sigma)\Rightarrow \max_{i,j\in \langle n\rangle}|h(i)-h(j)|<n.$$
By Proposition \ref{Proposition:EquiNumerous}
there exists unique $h_0\in \Motz_n(\theta,\sigma)$ 
such that $h_0\circ \theta-h_0=g$ and $h_0(1)=0$.
Let
$$\widehat{H}=\sum_{
\begin{subarray}{c}
h\in \Motz_n^N(\theta,\sigma)\\
\textup{\mbox{\scriptsize s.t.}}\,h'=g\;\;\;\;
\end{subarray}}\frac{1}{N}\left(\frac{h(1)}{N}\right)^{\nu}
=
\sum_{\begin{subarray}{c}
k\in \langle N\rangle\;\mbox{\scriptsize s.t.}\\
1\leq k+\min h_0\;\;\mbox{\scriptsize and}\\
k+\max h_0\leq N
\end{subarray}}\frac{1}{N}\left(\frac{k}{N}\right)^{\nu}
$$
where the second equality is justified by
\eqref{equation:MotzStability} and Proposition \ref{Proposition:EquiNumerous}.
Then we have
$$
|H-\widehat{H}|\leq \frac{n^2}{N}\;\;\mbox{and}\;\;
\left|-\widehat{H}+\sum_{k=1}^N \frac{1}{N}\left(\frac{k}{N}\right)^\nu\right|
\leq \frac{2n}{N}.
$$
Finally, we have evident inequalities
$$
\sum_{k=0}^{N-1}\frac{1}{N}\left(\frac{k}{N}\right)^{\nu}\leq
\int_0^1 t^{\nu}\,\dd t\leq \sum_{k=1}^{N}\frac{1}{N}\left(\frac{k}{N}\right)^{\nu}.
$$
Estimate \eqref{equation:FinalIntegrationBis} follows from the inequalities on the last two displayed lines. 
\qed

\section{Proof of Theorem \ref{Theorem:MainResult}}\label{section:TridiagonalCumulantCalc}

\subsection{Refinement of expansion \eqref{equation:MotzkinExpansionSteroid}}

\begin{Proposition}\label{Proposition:ProgressSoFar}
In the setup of Proposition \ref{Proposition:ToolCrush}  we have the yet more refined expansion
\begin{equation}\label{equation:ProgressSoFar}
\Mfrak_{\lambda,N}=
\sum_{(\sigma,A,h,\Gamma)}
\int2^{-|A|} \Ebold\left[(D^\Gamma f^h)(\zeta^h\star Q)\right]\, \PP^\Theta_\Gamma(\dd Q)
\end{equation}
where the sum is extended over quadruples $(\sigma,A,h,\Gamma)$
such that $\sigma\in \GJ_n(\theta)$, $A$ is a cycle-cutting of $\sigma$, $h\in \widetilde{\Motz}_n^N(\theta,\sigma)$
and $\Gamma=\LF(\sigma,A)$.
\end{Proposition}
\noindent The expansion \eqref{equation:ProgressSoFar} has interest beyond the scope of this paper.
Conceivably one could work out the $\frac{1}{N}$-expansion of the right side
and derive an alternate interpretation for the coefficients of the $\frac{1}{N}$-expansion
of $\Mfrak_{\lambda,N}$.

\proof By Proposition \ref{Proposition:ToolCrush}
we have
\begin{eqnarray*}
\Mfrak_{\lambda,N}&=&\sum_{\begin{subarray}{c}
(\Gamma,h)\in \Tree_n(\Theta)\times\Motz_n^N(\theta)\\
\mbox{\scriptsize s.t. \eqref{equation:MotzScrambleRedux1}---\eqref{equation:MotzScrambleRedux4} hold.}
\end{subarray}}\int \Ebold\left[(D^\Gamma f^h)(\zeta^h\star Q)\right]\, \PP^\Theta_\Gamma(\dd Q)\\
&=&\sum_{\Gamma\in \Tree^{\lf}_n(\Theta)}\sum_{\begin{subarray}{c}
h\in \Motz_n^N(\theta)\\
\mbox{\scriptsize s.t. \eqref{equation:MotzScrambleRedux1}---\eqref{equation:MotzScrambleRedux3} hold.}
\end{subarray}}\int \Ebold\left[(D^\Gamma f^h)(\zeta^h\star Q)\right]\, \PP^\Theta_\Gamma(\dd Q).
\end{eqnarray*}
By  Proposition \ref{Proposition:LoopCuttingGJ}
the formula \eqref{equation:ProgressSoFar} holds
with the sum is extended over quadruples 
$(\sigma,A,h,\Gamma)$ such that $\sigma\in \GJ_n(\theta)$,
$A$ is a cycle-cutting of $\sigma$, $h\in\Motz_n^N(\theta)$
satisfies \eqref{equation:MotzScrambleRedux1}---\eqref{equation:MotzScrambleRedux3}
and $\Gamma=\LF(\sigma,A)$.
Note that Lemma \ref{Lemma:BreakCycle} justifies the correction factor $2^{-|A|}$.
The formula \eqref{equation:ProgressSoFar} then holds as stated by 
Proposition \ref{Proposition:MotzScrambleInterp}.
\qed

\begin{Lemma}\label{Lemma:Cull}
Let $(\sigma,A,h,\Gamma)$ be a quadruple indexing a summand on the right side of \eqref{equation:ProgressSoFar}. Let
\begin{equation}\label{equation:SSAnotation}
S_0=J_h(0)\setminus \supp\,\sigma,\;\;
S_1=J_h(1)\setminus \supp\,\sigma,\;\;
A_1=J_h(1)\cap A.
\end{equation}
Then we have
\begin{equation}
\label{equation:CombinatorialDifferentiation}
D^\Gamma f^h
=\left(\prod_{i\in S_0}z_{i0}\right)
\left(\prod_{i\in S_1}\sum_{j=1}^{2h(i)}\frac{z_{ij}^2}{2}
\right)
\left(\prod_{i\in A_1}
\sum_{j=1}^{2 h(i)}z_{ij}z_{\sigma(i),j}\right).
\end{equation}
\end{Lemma}

\proof  Let 
$$\sigma=\sigma_1\cdots \sigma_p\tau_1\cdots \tau_q$$ be the canonical factorization
of $\sigma$ into disjoint cycles, with the factors sorted so that
$$\bigcup_{\alpha=1}^p\supp\, \sigma_\alpha\subset J_h(1)\;\;\mbox{and}\;\;\bigcup_{\beta=1}^q\supp\, \tau_\beta\subset J_h(0).$$
Such a sorting is possible because $h\circ \sigma=h$.
Note that each permutation $\tau_\beta$ is necessarily a transposition
since $i\in J_h(0)\Rightarrow\sigma^2(i)=i$.
For $\alpha=1,\dots,p$ and $\beta=1,\dots,q$ let
$$\{a_\alpha\}=A\cap \supp \,\sigma_\alpha,\;M_\alpha=2h(a_\alpha),\;\Gamma_\alpha=\LF(\sigma_\alpha,\{a_\alpha\})\;\;\mbox{and}\;\;e_\beta=\supp\,\tau_\beta.$$
Then
$$\Gamma=\bigcup_{\alpha=1}^p\Gamma_\alpha\cup \bigcup_{\beta=1}^q\{e_\beta\}$$
is the decomposition of $\Gamma$ into connected components.
Note that since $h$ is constant on $\sigma$-orbits, $M_\alpha$ is the value of $2h$ on $\supp\,\sigma_\alpha$.
Then have a factorization
\begin{eqnarray*}
D^\Gamma f^h&=&\left(\prod_{i\in S_0}z_{i0}\right)
\left(\prod_{i\in S_1}\sum_{j=1}^{2h(i)}\frac{z_{ij}^2}{2}\right)\\
&&\times \left(\prod_{\alpha=1}^p D^{\Gamma_\alpha}\prod_{i\in \supp\,\sigma_\alpha}\sum_{j=1}^{M_\alpha}\frac{z_{ij}^2}{2}\right)\left(\prod_{\beta=1}^q D_{e_\beta}\prod_{i\in e_\beta}z_{i0}\right).
\end{eqnarray*}
It is easy to see that 
$$D_{e_\beta}\prod_{i\in e_\beta}z_{i0}=1.$$
To finish the proof we need only evaluate
$$D^{\Gamma_\alpha}\prod_{i\in \supp\,\sigma_\alpha}\sum_{j=1}^{M_\alpha}\frac{z_{ij}^2}{2}.$$
For the latter purpose we note the formula
$$
\left(\sum_{j=0}^{2N} \frac{\partial^2}{\partial z_{i_2j}\partial z_{i_3j}}\right)
\left[
\left(\frac{|\{i_1,i_2\}|}{2}\sum_{j=1}^{M} z_{i_1j}z_{i_2j}\right)\left(
\sum_{j=1}^{M}\frac{z_{i_3j}^2}{2}\right)\right]
=\sum_{j=1}^{M}z_{i_1j}z_{i_3j}
$$
holding for $i_1,i_2,i_3\in \langle n\rangle$ such that $i_3\not\in \{i_1,i_2\}$
and $1\leq M\leq 2N$.
Using this relation and induction one can finish the proof. We omit the remaining details. \qed

\subsection{Application of the Marcinkiewicz-Zygmund inequality}
\subsubsection{The Marcinkiewicz-Zygmund inequality}
For a real random variable $Z$ and $p\in [1,\infty)$,
let $\norm{Z}_p=(\Ebold|Z|^p)^{1/p}$. Now fix $p\in[1,\infty)$
and let $X_1,\dots,X_n$ be independent real random variables with finite $L^p$-norms,
each of mean zero. Then the {\em Marcinkiewicz-Zygmund inequality}
is the assertion that
$$
A_p\norm{\left(\sum_{i=1}^N X_i^2\right)^{1/2}}_p
\leq \norm{\sum_{i=1}^N X_i}_p\leq B_p\norm{\left(\sum_{i=1}^N X_i^2\right)^{1/2}}_p
$$
for positive constants $A_p$ and $B_p$ depending only on $p$.
See \cite[p. 386]{ChowTeich} for a textbook treatment of this inequality.
For $p\geq 2$ via the Minkowski inequality we deduce the relatively crude inequality
$$\norm{\sum_{i=1}^N X_i}_p\leq \frac{K_p}{2}
\left(\sum_{i=1}^N \norm{X_i}_p^2\right)^{1/2}\leq \frac{K_p}{2}\sqrt{N}\max_{i=1}^N \norm{X_i}_p
$$
for a constant $K_p$ depending only on $p$.
Let  $T_1,\dots,T_N$ be independent real random variables with finite $L^p$-norms
which might not all be of mean zero. Finally we have a bound
\begin{equation}\label{equation:MarcinkiewiczZygmundCrude}
\norm{\sum_{i=1}^NT_i-\Ebold\sum_{i=1}^N T_i}_{p}
\leq
\frac{K_p}{2}\sqrt{N}\max_{i=1}^N\norm{T_i-\Ebold T_i}_{p}
\leq K_p \sqrt{N}\max_{i=1}^N\norm{T_i}_{p}
\end{equation}
which is all we need going forward.

\begin{Lemma}\label{Lemma:KillZero}
Let $(\sigma,A,h,\Gamma)$ be a quadruple indexing a summand on the right side of \eqref{equation:ProgressSoFar},
and let $S_0$, $S_1$ and $A_1$ be as defined on line \eqref{equation:SSAnotation}.
 Let $\ell=\ell(\theta)$. 
Fix $Q\in \Qfrak_n$ arbitrarily.
We have
\begin{equation}\label{equation:NuggetBis}
\left|\Ebold[(D^\Gamma f^h)(\zeta^h\star Q)]-\indicator{S_0=\emptyset}
\cdot\prod_{i\in A_1}2Q(i,\sigma(i))
\cdot \prod_{i\in S_1\cup A_1}h(i)\right|\leq c  N^{\frac{n}{2}-\ell+\frac{1}{2}}
\end{equation}
for a constant $c$ depending only on $n$. 
\end{Lemma}
\proof
Let $S=S_0\cup S_1\cup A_1$.  For $i\in S$ let
\begin{equation}
Z_i=\left\{\begin{array}{rl}
z_{i0}\bigg\vert_{\zeta^h\star Q}&\mbox{if $i\in S_0$,}\\
\displaystyle\sum_{j=1}^{2h(i)}\frac{z_{ij}^2}{2}\bigg\vert_{\zeta^h\star Q}&\mbox{if $i\in S_1$,}\\
\displaystyle\sum_{j=1}^{2h(i)}z_{ij}z_{\sigma(i),j}\bigg\vert_{\zeta^h\star Q}&\mbox{if $i\in A_1$.}
\end{array}\right.
\end{equation}
By the definitions and Lemma \ref{Lemma:Cull} we have
\begin{equation}
\Ebold
\prod_{i\in S}Z_i=
\Ebold [(D^\Gamma f)(\zeta^h\star Q)].
\end{equation}
We also have
\begin{equation}
\Ebold Z_i=\left\{\begin{array}{rl}
0&\mbox{if $i\in S_0$,}\\
h(i)&\mbox{if $i\in S_1$,}\\
2Q(i,\sigma(i))h(i)&\mbox{if $i\in A_1$}
\end{array}\right.
\end{equation}
by using the fact that by definition $\zeta^h\star Q$ is a centered
Gaussian random vector with covariances 
$$\Ebold (\zeta\star Q)_{ij}(\zeta\star Q)_{i'j'}=\delta_{h(i),h(i')}\delta_{jj'}Q(i,i').$$
Using this same covariance information,
 the general bound \eqref{equation:MarcinkiewiczZygmundCrude} recalled above 
and the fact that $\zeta^h$ is a Gaussian random vector,
we also have
\begin{equation}
\norm{Z_i-\Ebold Z_i}_{2n}\leq \gamma\left\{\begin{array}{rl}
1&\mbox{if $i\in S_0$,}\\
\sqrt{N}&\mbox{if $i\in S_1\cup A_1$,}
\end{array}\right.
\end{equation}
where the constant $\gamma \geq 1$ depends only on  $n$.
Finally we have
$$\norm{\prod_{i\in S}Z_i-\prod_{i\in S} \Ebold Z_i}_2
\leq\sum_{i\in S}\norm{Z_i-\Ebold Z_i}_{2n}\prod_{i'\in S\setminus \{i\}}\norm{Z
_{i'}}_{2n}\leq cN^{\min(|S\setminus S_0|,|S|-\frac{1}{2})}$$
where $c$ depends only on $n$,
whence estimate \eqref{equation:NuggetBis} by Lemma \ref{Lemma:BalanceBooks}. \qed

\subsection{Closing arguments to prove Theorem \ref{Theorem:MainResult}}
From \eqref{equation:NuggetBis}, by integrating on both sides against $\PP_\Gamma^\Theta$
and using Jensen's inequality, along with formula \eqref{equation:QInt} from Proposition \ref{Proposition:LastCut} and formula \eqref{equation:SpinUp}, we deduce
the inequality
\begin{equation}\label{equation:Nugget}
\left|2^{-|A|}\int\,\Ebold[(D^\Gamma f^h)(\zeta^h\star Q)]\dd \PP_\Gamma^\Theta(Q)-
\frac{\indicator{S_0=\emptyset}}{\mbold(\sigma)}
\prod_{i\in S_1\cup A_1}h(i)
\right|\leq c_1  N^{\frac{n}{2}-\ell+\frac{1}{2}}
\end{equation}
where the constant $c_1$ depends only on $n$.
After using \eqref{equation:CrudeMotzBound} and \eqref{equation:Nugget}
to approximate the right side of \eqref{equation:ProgressSoFar}, we obtain the approximation
$$
\left|\Mfrak_{\lambda,N}-
\sum_{(\sigma,A,h,\Gamma)}
\frac{\indicator{J_h(0)\subset \supp\,\sigma}}{\mbold(\sigma)}
\prod_{i\in J_h(1)\setminus((\supp\,\sigma)\setminus A)}h(i)\right|\leq c_2N^{\frac{n}{2}-\ell+\frac{3}{2}}
$$
where $c_2$ depends only on $n$ and the sum is extended over the same family of quadruples $(\sigma,A,h,\Gamma)$ as in \eqref{equation:ProgressSoFar}.
Using \eqref{equation:UltraTrivial}, \eqref{equation:FinalIntegration}, and again \eqref{equation:CrudeMotzBound} we then get a further
approximation
$$\left|\Mfrak_{\lambda,N}-\frac{N^{\frac{n}{2}-\ell+2}}{\frac{n}{2}-\ell+2}
\cdot\sum_{\sigma\in \GJ_n(\theta)}\sum_{\begin{subarray}{c}
\mbox{\scriptsize cycle-cuttings}\\
\mbox{\scriptsize $A$ of $\sigma$}
\end{subarray}}\frac{|\dMotz_n(\theta,\sigma)|}{\mbold(\sigma)}\right|\leq c_3N^{\frac{n}{2}-\ell+\frac{3}{2}}$$
where $c_3$ depends only on $n$. Note that the inner sum over cycle-cuttings $A$
of $\sigma$ is  in effect canceled by the factor $1/\mbold(\sigma)$.
Thus the last estimate in conjunction with limit formula \eqref{equation:tHooft} 
and the definition of $\GJdM_n(\theta)$ proves Theorem \ref{Theorem:MainResult}. \qed

\end{document}